\documentclass{article}
\usepackage{amsmath}
\usepackage{amssymb}
\begin{document}
\def\e#1\e{\begin{equation}#1\end{equation}}
\def\ea#1\ea{\begin{align}#1\end{align}}
\def\eq#1{{\rm(\ref{#1})}}
\newtheorem{thm}{Theorem}[section]
\newtheorem{lem}[thm]{Lemma}
\newtheorem{prop}[thm]{Proposition}
\newtheorem{cor}[thm]{Corollary}
\newenvironment{result}[1]{\medskip\noindent{\bf #1}\ \it}{\rm\medskip}
\newenvironment{dfn}{\medskip\refstepcounter{thm}
\noindent{\bf Definition \thesection.\arabic{thm}\ }}{\medskip}
\newenvironment{proof}[1][,]{\medskip\ifcat,#1
\noindent{\it Proof.\ }\else\noindent{\it Proof of #1.\ }\fi}
{\hfill$\square$\medskip}
\def\dim{\mathop{\rm dim}}
\def\Re{\mathop{\rm Re}}
\def\Im{\mathop{\rm Im}}
\def\supp{\mathop{\rm supp}}
\def\vol{\mathop{\rm vol}}
\def\ge{\geqslant} 
\def\le{\leqslant} 
\def\R{\mathbin{\mathbb R}}
\def\ha{\textstyle\frac{1}{2}}
\def\cal{\mathcal}
\def\al{\alpha}
\def\be{\beta}
\def\ga{\gamma}
\def\de{\delta}
\def\la{\lambda}
\def\ep{\epsilon}
\def\ze{\zeta}
\def\si{\sigma}
\def\De{\Delta}
\def\La{\Lambda}
\def\Ga{\Gamma}
\def\Si{\Sigma}
\def\om{\omega}
\def\d{{\rm d}}
\def\ts{\textstyle}
\def\sst{\scriptscriptstyle}
\def\sm{\setminus}
\def\iy{\infty}
\def\ra{\rightarrow}
\def\hookra{\hookrightarrow}
\def\t{\times}
\def\na{\nabla}
\def\ti{\tilde}
\def\ov{\overline}
\def\ovB{\,\overline{\!B}}
\def\bigti{\raise-4pt\hbox{$\widetilde{\hskip 11pt}$}}
\def\ms#1{\vert#1\vert^2}
\def\bms#1{\bigl\vert#1\bigr\vert^2}
\def\md#1{\vert #1 \vert}
\def\bmd#1{\big\vert #1 \big\vert}
\def\nm#1{\Vert #1 \Vert}
\def\an#1{\langle#1\rangle}
\def\cnm#1#2{\Vert #1 \Vert_{C^{#2}}} 
\def\lnm#1#2{\Vert #1 \Vert_{L^{#2}}} 
\def\snm#1#2#3{\Vert #1 \Vert_{L^{#2}_{#3}}} 
\def\bnm#1{\bigl\Vert #1 \bigr\Vert}
\def\bcnm#1#2{\bigl\Vert #1 \bigr\Vert_{C^{#2}}} 
\def\blnm#1#2{\bigl\Vert #1 \bigr\Vert_{L^{#2}}} 
\def\bsnm#1#2#3{\bigl\Vert #1 \bigr\Vert_{L^{#2}_{#3}}} 
\title{Constant Scalar Curvature Metrics \\ on Connected Sums}
\author{Dominic Joyce \\ Lincoln College, Oxford, OX1 3DR, U.K. \\
{\tt dominic.joyce@lincoln.ox.ac.uk}}
\date{}
\maketitle
\begin{abstract}
The {\it Yamabe problem} (proved in 1984) guarantees the
existence of a metric of constant scalar curvature in each 
conformal class of Riemannian metrics on a compact manifold
of dimension $n\ge 3$, which minimizes the total scalar 
curvature on this conformal class.

Let $(M',g')$ and $(M'',g'')$ be compact Riemannian $n$-manifolds. 
We form their {\it connected sum} $M'\# M''$ by removing small balls 
of radius $\epsilon$ from $M'$, $M''$ and gluing together the 
${\mathcal S}^{n-1}$ boundaries, and make a metric $g$ on $M'\# M''$ 
by joining together $g',g''$ with a partition of unity.

In this paper we use analysis to study metrics with constant scalar 
curvature on $M'\# M''$ in the conformal class of $g$. By the Yamabe 
problem, we may rescale $g'$ and $g''$ to have constant scalar 
curvature 1, 0 or $-1$. Thus there are 9 cases, which we handle separately. 

We show that the constant scalar curvature metrics either develop
small `necks' separating $M'$ and $M''$, or one of $M',M''$ is 
crushed small by the conformal factor. When both sides have positive 
scalar curvature we find three metrics with scalar curvature 1 in the 
same conformal class.


\noindent{\bf 2000 MSc classifications:} primary 58E11, secondary 53C20.
\end{abstract}
\section{Introduction}
\label{s1}

Let $(M',g')$ and $(M'',g'')$ be compact manifolds of dimension
$n\ge 3$. The {\it connected sum} $M=M'\# M''$ is the result of
removing a small ball $B^n$ from each manifold, and joining the
resulting manifolds at their common boundary ${\cal S}^{n-1}$.
By smoothly joining $g'$ and $g''$ we can also construct a
1-parameter family of metrics $g_t$ on $M$ for $t\in(0,\de)$,
where $t$ measures the radius of the excised balls~$B^n$.

In this paper we suppose $g'$ and $g''$ have constant scalar
curvature (not necessarily the same value), and study (some of)
the metrics $\ti g_t$ with constant scalar curvature in the
conformal class of $g_t$ for small $t$. Our method is to write
down explicit metrics whose scalar curvature is close to constant,
and show using analysis that they can be adjusted by a small
conformal change to give metrics with constant scalar curvature.

By the proof of the Yamabe problem, every conformal class on $M$
contains a metric with constant scalar curvature. Our proofs are
simpler than the solution of the Yamabe problem, as our problem
is much easier, and they have the advantage of giving a good grasp
of what the Yamabe metrics actually look like as the underlying
conformal manifold decays into a connected sum. We can for instance
say that one obvious sort of behaviour, that of developing long
`tubes' resembling ${\cal S}^{n-1}\t\R$ does not happen, but that
small `pinched necks' may develop instead, or else $M'$ or $M''$
may be crushed very small by the conformal factor.

Other authors have also considered the scalar curvature of metrics on
connected sums. Two important early papers on manifolds with positive
scalar curvature are Gromov and Lawson \cite{GrLa} and Schoen and Yau
\cite{ScYa}, and both show \cite[Th.~A]{GrLa}, \cite[Cor.~3]{ScYa}
that the connected sum of two manifolds with positive scalar curvature
carries metrics with positive scalar curvature.

Writing after the solution of the Yamabe problem, Kobayashi \cite{Koba}
defines the {\it Yamabe number} of a manifold to be the supremum over
conformal classes of the {\it Yamabe invariant} defined in \S\ref{s12},
and proves an inequality \cite[Th.~2(a)]{Koba} between the Yamabe number
of two manifolds and the Yamabe number of their connected sum. His
formula can be related to our results below in the limit~$t\ra 0$.

The first version of this work formed part of the author's D.Phil.
thesis \cite{Joyc} in 1992, supervised by Simon Donaldson. In the
interval between then and the publication of the present paper
several papers on similar topics have appeared, and we mention in
particular Mazzeo, Pollack and Uhlenbeck \cite{MPU}. They study
connected sums of compact or noncompact manifolds with constant
positive scalar curvature, and the part of Theorem \ref{t331}
below dealing with the metrics of \S\ref{s22} also follows from
their main result.

The remainder of this section goes over some necessary backgound material.
In \S\ref{s2} we define two families of metrics $g_t$ upon connected sums,
make estimates of the scalar curvature of these metrics, and prove a uniform
bound on a Sobolev constant for a particular embedding of Sobolev spaces.

The main existence results for metrics of constant scalar curvature
conformal to the metrics $g_t$ above are proved in \S\ref{s3}. We
begin with a quite general existence proof using a sequence method
and then apply it, first to the case of scalar curvature $-1$, and
then to the case of scalar curvature 1. The former is simple, but
the latter is more difficult, and the proof of a result on the
eigenvalues of the Laplacian $\De$ on the metrics $g_t$ in the
positive case has been deferred until the appendix.

In \S\ref{s4} we deal with the cases left over from \S\ref{s3}, which are
the connected sums involving manifolds of zero scalar curvature. To do so
requires some additions to the methods of \S\ref{s3}, as the problem of
rescaling a metric with scalar curvature close to $\pm 1$ to get exactly
$\pm 1$ is different to the problem of rescaling a metric with scalar
curvature close to zero to get scalar curvature a small but unknown
constant. We shall see that each combination of positive, negative and
zero scalar curvature manifolds has distinctive features.
\medskip

\noindent{\it Acknowledgements.} The author would like to thank
Simon Donaldson for many useful conversations, and Rafe Mazzeo
and the referees for helpful comments.

\subsection{Analytic preliminaries}
\label{s11}

We define Lebesgue and Sobolev spaces, principally to establish notation.
An introduction to them may be found in Aubin \cite[\S 2]{Aubi}. Let
$(M,g)$ be a Riemannian manifold. For $q\ge 1$, the {\it Lebesgue space}
$L^q(M)$ is the set of locally integrable functions $u$ on $M$ for which
the norm
\begin{equation*}
{}^g\lnm{u}{q}=\lnm{u}{q}=\left(\int_M\md{u}^q\,\d V_g\right)^{1/q}
\end{equation*}
is finite. Here $\d V_g$ is the volume form of the metric $g$.
Then $L^q(M)$ is a Banach space (with the convention that two
functions are equal if they differ only on a null set) and
$L^q(M)$ a Hilbert space. Let $r,s,t\ge 1$ with $1/r=1/s+1/t$.
If $\phi\in L^s(M)$ and $\psi\in L^t(M)$ then $\phi\psi\in L^r(M)$
and $\lnm{\phi\psi}{r}\le\lnm{\phi}{s}\lnm{\psi}{t}$; this is
{\it H\"older's inequality}.

Let $q\ge 1$ and let $k$ be a nonnegative integer. The {\it Sobolev space}
$L^q_k(M)$ is the subspace of $u\in L^q(M)$ such that $u$ is $k$ times
weakly differentiable and $\md{\na^iu}\in L^q(M)$ for $i\le k$. The
{\it Sobolev norm} on $L^q_k(M)$ is
\begin{equation*}
{}^g\snm{u}{q}{k}=\snm{u}{q}{k}=
\left(\sum_{i=0}^k\int_M\md{\na^iu}^q\,\d V_g\right)^{1/q}.
\end{equation*}
This makes $L^q_k(M)$ into a Banach space, and $L^2_k(M)$ a Hilbert
space. Since we will often have to consider different metrics on the
same manifold, we will use the superscript notation ${}^g\lnm{\,.\,}{q}$
and ${}^g\snm{\,.\,}{q}{k}$ to mean Lebesgue and Sobolev norms computed
using the metric~$g$.

We also write $C^k(M)$ for the vector space of continuous, bounded
functions on $M$ with $k$ continuous, bounded derivatives, and
$C^\iy(M)=\bigcap_{k\ge 0}C^k(M)$. The {\it Sobolev Embedding
Theorem} \cite[Th.~2.20]{Aubi} says that if $\frac{1}{r}\ge\frac{1}{q}-
\frac{k}{n}$ then $L^q_k(M)$ is continuously embedded in $L^r(M)$
by inclusion, and if $-\frac{r}{n}>\frac{1}{q}-\frac{k}{n}$ then
$L^q_k(M)$ is continuously embedded in $C^r(M)$ by inclusion.

\subsection{The Yamabe problem}
\label{s12}

We now briefly discuss the {\it Yamabe problem}, concerning the
existence of metrics of constant scalar curvature in a conformal
class on a compact manifold. An introduction to the problem and
its solution may be found in the survey paper \cite{LePa} by Lee
and Parker.

Fix a dimension $n\ge 3$. Throughout the paper we shall use the
notation
\e
a=\frac{4(n-1)}{n-2},\quad b=\frac{4}{n-2}\quad\text{and}\quad
p=\frac{2n}{n-2}
\label{s1eq1}
\e
Let $M$ be a compact manifold of dimension $n$. Define a
functional $Q$ upon the set of Riemannian metrics on $M$ by
\begin{equation*}
Q(g)=\frac{\int_MS\,\d V_g}{\vol(M)^{2/p}},
\end{equation*}
where $S$ is the scalar curvature of $g$. Then $Q$ is known as
the {\it total scalar curvature}, or {\it Hilbert action}.

Let $g,\ti g$ be metrics on $M$. We say that $g,\ti g$ are
{\it conformal} if $\ti g=\phi g$ for some smooth {\it conformal
factor} $\phi:M\ra(0,\iy)$. Usually we write $\phi=\psi^{p-2}$ for
smooth $\psi:M\ra(0,\iy)$. The {\it conformal class} $[g]$ of $g$
is the set of metrics $\ti g$ on $M$ conformal to $g$. Applying the
calculus of variations to the restriction of $Q$ to a conformal
class, one finds that a metric $g$ is a stationary point of $Q$
on $[g]$ if and only if it has {\it constant scalar curvature}.

Let the scalar curvatures of $g,\ti g$ be $S,\ti S$, and write
$\ti g=\psi^{p-2}g$ for smooth $\psi:M\ra(0,\iy)$. Then
\e
\ti S=\psi^{1-p}\left(a\De\psi+S\psi\right).
\label{s1eq2}
\e
Thus $\ti g=\psi^{p-2}g$ has constant scalar curvature $\nu$ if
and only if $\psi$ satisfies the {\it Yamabe equation}:
\e
a\De\psi+S\psi=\nu\md{\psi}^{p-1}.
\label{s1eq3}
\e

One can show using H\"older's inequality that $Q$ is bounded
below on a conformal class. So we define
\begin{equation*}
\la\bigl([g]\bigr)=\inf\bigl\{Q(\ti g):\ti g\text{\ is
conformal to $g$}\bigr\}.
\end{equation*}
This constant $\la\bigl([g]\bigr)$ is an invariant of the conformal
class $[g]$ of $g$ on $M$, called the {\it Yamabe invariant}.
\medskip

\noindent{\bf The Yamabe problem \cite{Yama}.} {\it Given a compact 
Riemannian $n$-manifold\/ $(M,g)$, find a metric
$\ti g$ conformal to $g$ which minimizes $Q$ on $[g]$, so that
$Q(\ti g)=\la\bigl([g]\bigr)$. Then $\ti g$ is a stationary point
of\/ $Q$ on $[g]$, and so has constant scalar curvature.}
\medskip

A solution $\ti g$ of the Yamabe problem is called a {\it Yamabe
metric}. The problem was posed in 1960 by Yamabe \cite{Yama}, who
gave a proof that such a metric $\ti g$ always exists. His idea was
to choose a minimizing sequence for $Q$ on $[g]$, and show that a
subsequence converges to a smooth minimizer $\ti g$ for $Q$.
Unfortunately the proof contained an error, discovered by Trudinger
\cite{Trud}. The proof was eventually repaired by Trudinger, Aubin,
Schoen and Yau.

If $\la\bigl([g]\bigr)\le 0$ then constant scalar curvature
metrics in $[g]$ are unique up to homothety \cite[p.~135]{Aubi},
and are Yamabe metrics. Thus, if $\la\bigl([g]\bigr)<0$ then $[g]$
contains a unique metric $\ti g$ with scalar curvature $-1$, and if
$\la\bigl([g]\bigr)=0$ then $[g]$ contains a unique metric $\ti g$
with scalar curvature 0 and volume 1. However, if $\la\bigl([g]\bigr)>0$
then constant scalar curvature metrics in $[g]$ are not necessarily
unique up to homothety, and may be higher stationary points of $Q$
rather than Yamabe metrics. Thus, $[g]$ may contain several different
metrics with scalar curvature 1. We will see an example of this
in~\S\ref{s33}.

In \S\ref{s3} and \S\ref{s4}, weak solutions in $L^2_1(M)$ to the
Yamabe equation \eq{s1eq3} will be constructed for certain special
compact manifolds $M$. For these solutions to give metrics of
constant scalar curvature, it is necessary that they be not just
$L^2_1$ solutions, but $C^\iy$ solutions. We quote a result of
Trudinger \cite[Th.~3, p.~271]{Trud} showing this is the case.

\begin{prop} Let\/ $(M,g)$ be a compact Riemannian $n$-manifold
with scalar curvature $S$, and\/ $u$ a weak $L^2_1$ solution of\/
$a\De u+Su=\ti S\md{u}^{(n+2)/(n-2)}$, for $\ti S\in C^\iy(M)$.
Then $u\in C^2(M)$, and where $u$ is nonzero it is~$C^\iy$.
\label{p122}
\end{prop}

\subsection{Stereographic projections}
\label{s13}

Let $(M,g)$ be a compact Riemannian manifold with positive scalar
curvature. Then by \cite[Lem.~6.1, p.~63]{LePa} the operator
$a\De+S$ in \eq{s1eq3} admits a {\it Green's function} $\Ga_m$
for each $m\in M$, which is unique and strictly positive. That
is, $\Ga_m:M\sm\{m\}\ra(0,\iy)$ is smooth with $(a\De+S)\Ga_m=0$
on $M\sm\{m\}$, and $\Ga_m$ is singular at $m$ with $(a\De+S)
\Ga_m=\de_m$ on $M$ in the sense of distributions.

\begin{dfn} Let $(M,g)$ be a compact Riemannian manifold with
positive scalar curvature. For $m\in M$ define the metric $\hat g
=\Ga_m^{p-2}g$ on $\hat M=M\sm\{m\}$, where $\Ga_m$ is the Green's
function of $a\De+S$ at $m$. Then $\hat g$ has {\it zero scalar
curvature} by \eq{s1eq2}. We call $(\hat M,\hat g)$ the
{\it stereographic projection} of $M$ from~$m$.
\label{d132}
\end{dfn}

The stereographic projection of ${\cal S}^n$ is $\R^n$. In
general, stereographic projections $(\hat M,\hat g)$ are
{\it asymptotically flat}, that is, their noncompact ends
resemble $\R^n$ in a way made precise in \cite[Def.~6.3]{LePa}.
Following Lee and Parker \cite[p.~64]{LePa}, we use the notation
that $\phi=O'(\md{x}^q)$ means $\phi=O(\md{x}^q)$ and $\na\phi=
O(\md{x}^{q-1})$, and $\phi=O'(\md{x}^q)$ means $\phi=O(\md{x}^q)$,
$\na\phi=O(\md{x}^{q-1})$ and $\na^2\phi=O(\md{x}^{q-2})$, for small
or large $x$ in $\R^n$, depending on the context. For simplicity we
shall suppose that $g$ is {\it conformally flat} near $m$ in $M$.
In this case the asymptotic expansion of the metric $\hat g$ is
particularly simple,~\cite[Th.~6.5(a)]{LePa}.

\begin{prop} Let\/ $(M,g)$ be a compact Riemannian manifold with\/
positive scalar curvature, and suppose $m\in M$ has a neighbourhood
that is conformally flat. Let\/ $(\hat M,\hat g)$ be the 
stereographic projection of\/ $M$ from $m$. Then there exist\/
$\mu\in\R$ and asymptotic coordinates $\{x^i\}$ on $\hat M$
with respect to which\/ $\hat g_{ij}$ satisfies
\begin{equation*}
\hat g_{ij}(x)=\psi^{p-2}(x)\de_{ij},\;\>\text{where}\;\>
\psi(v)=1+\mu\md{x}^{2-n}+O''(\md{x}^{1-n})\;\>\text{for large $\md{x}$.}
\end{equation*}
\label{p134}
\end{prop}

This constant $\mu$ is proportional to an important invariant of
asymptotically flat manifolds called the {\it mass}, which is
defined and studied by Bartnik \cite{Bart}. For physical reasons
it was conjectured that the mass must be nonnegative for complete
asymptotically flat manifolds of nonnegative scalar curvature, and
zero only for flat space. This was proved for {\it spin manifolds}
by Witten \cite{Witt}, whose proof was generalized to $n$ dimensions
by Bartnik~\cite[\S 6]{Bart}.

A proof of the general case when $n\le 7$ is given by Schoen
\cite[\S 4]{Scho2}. He also claims the result for all dimensions
\cite[p.~481]{Scho1}, \cite[p.~145]{Scho2}, but the proof has not
been published. This `positive mass theorem' in dimensions 3,4 and 5
was important in the completion of the proof of the Yamabe problem by
Schoen in 1984 \cite{Scho1}. We will only need the case when the metric
$g$ of $M$ is conformally flat about $m$, which we state as follows:

\begin{thm} In the situation of Proposition \ref{p134}, (perhaps with
the additional assumption that\/ $M$ is spin or $n\le 7$), we have
$\mu\ge 0$, with equality if and only if\/ $M$ is ${\cal S}^n$, and
its metric is conformal to the round metric.
\label{t135}
\end{thm}

\section{Glued metrics on connected sums}
\label{s2}

In this section, we shall define two families of metrics $g_t$ on the 
connected sum $M$ of constant scalar curvature Riemannian manifolds 
$(M',g')$ and $(M'',g'')$. The first family, in \S\ref{s21}, is made
by choosing $M'$ of constant scalar curvature $\nu$, cutting out a
small ball, and gluing in a stereographic projection of $M''$,
homothetically shrunk very small. The second family, in \S\ref{s22},
is made by choosing $M'$ and $M''$ both of constant scalar curvature
$\nu$, and joining them by a small `neck' of zero scalar curvature.
The relation between these families when $\nu=1$ is discussed
in~\S\ref{s23}.

We finish with two results, Propositions \ref{p231} and \ref{p241},
about the families of metrics. The first gives explicit bounds for
their scalar curvature, to determine how good an approximation to
constant scalar curvature they are. The second shows that the
Sobolev constant for a certain Sobolev embedding of function
spaces can be given a bound independent of the width of the neck,
for small values of this parameter.

\subsection{Combining constant and positive scalar curvature}
\label{s21}

Let $(M',g')$ be a compact Riemannian $n$-manifold with constant
scalar curvature $\nu$. Applying a homothety to $g'$ if necessary,
we may assume that $\nu=1,0$ or $-1$. For simplicity, assume that
$m'\in M'$ has a neighbourhood in which $g'$ is conformally flat;
this assumption will be dropped in \S\ref{s34}. Then $M'$ contains
a ball $B'$ about $m'$, with a diffeomorphism $\Phi'$ from
$B_{\de}(0)\subset\R^n$ to $B'$ for some $\de\in(0,1)$, such that
$\Phi'(0)=m'$ and $(\Phi')^*(g')=(\psi')^{p-2}h$ for some function
$\psi'$ on $B_{\de}(0)$, where $h$ is the standard metric on $\R^n$.
By choosing a different conformal identification with $B_{\de}(0)$
if necessary, we may suppose $\psi'(0)=1$ and $\d\psi'(0)=0$,
so that $\psi'(v)=1+O'(\ms{v})$ in the notation of \S\ref{s13}.

Let $(M'',g'')$ be a compact Riemannian $n$-manifold with
scalar curvature 1. As with $M'$, suppose $m''\in M''$ has
a neighbourhood in which $g''$ is conformally flat. Let
$(\hat M,\hat g)$ be the stereographic projection of $M''$ from
$m''$, as in Definition \ref{d132}, so that $\hat M=M''\sm\{m''\}$,
and $\hat g$ is asymptotically flat with zero scalar curvature,
and is conformal to~$g''$.

By Proposition \ref{p134}, there is an immersion $\Xi'':\R^n\sm
\ovB_R(0)\ra\hat M$ for some $R>0$, whose image is the complement
of a compact set in $\hat M$, such that $(\Xi'')^*(\hat g)=\xi^{p-2}h$,
where $\xi$ is a smooth function on $\R^n\sm\ovB_R(0)$ satisfying
$\xi(v)=1+O'(\md{v}^{2-n})$. By making $\de$ smaller if necessary,
we set~$R=\de^{-4}$.

A family of metrics $\{g_t:t\in(0,\de)\}$ on $M=M'\# M''$ will now
be written down. For any $t\in(0,\de)$, define $M$ and the conformal
class of $g_t$ by
\begin{equation*}
M=\Bigl(M'\sm\Phi'\bigl[\ovB_{t^2}(0)\bigr]\Bigr)\amalg 
\Bigl(\hat M\sm\bigl(\Xi''\bigl[\R^n\sm
B_{t^{-5}}(0)\bigr]\bigr)\Bigr)\Bigl/\bigti_t,
\end{equation*} 
where $\bigti_t$ is the equivalence relation defined by
\begin{equation*}
\Phi'[v]\,\,\bigti_t\,\,\Xi''[t^{-6}v]
\quad\text{whenever $v\in\R^n$ and $t^2<\md{v}<t$.}
\end{equation*}

The conformal class $[g_t]$ of $g_t$ is then given by the restriction 
of the conformal classes of $g'$ and $\hat g$ to the open sets of $M',M''$ 
that make up $M$; this definition makes sense because the conformal 
classes agree on the annulus of overlap where the two open sets are 
glued by $\bigti_t$. Let $A_t$ be this annulus in $M$. Then $A_t$ is 
diffeomorphic via $\Phi'$ to the annulus $\{v\in\R^n:t^2<\md{v}<t\}$ 
in~$\R^n$.

To define a metric $g_t$ within the conformal class just given,
take $g_t=g'$ on the component of $M\sm A_t$ coming from $M'$,
and $g_t=t^{12}\hat g$ on the component coming from $\hat M$.
It remains to choose a conformal factor on $A_t$, which is
identified with $\{v\in\R^n:t^2<\md{v}<t\}$ by $\Phi'$. Here
we set $(\Phi')^*(g_t)=\psi_t^{p-2}h$, where $\psi_t$ is
defined below.

The conditions for smoothness at the edges of the annulus $A_t$ are
that $\psi_t(v)$ should join smoothly onto $\psi'(v)$ at $\md{v}=t$,
and onto $\xi(t^{-6}v)$ at $\md{v}=t^2$. Choose a $C^\iy$ function
$\si:\R\ra[0,1]$, that is 0 for $x\ge 2$ and 1 for $x\le 1$ and
strictly decreasing in $[1,2]$. Let $\be_1(v)=\si(\log\md{v}/\log t)$
and $\be_2(v)=1-\be_1(v)$ for all $v\in\R^n$ with $t^2<\md{v}<t$.
Finally, define $\psi_t$ by
\e
\psi_t(v)=\be_1(v)\psi'(v)+\be_2(v)\xi(t^{-6}v).
\label{s2eq1}
\e

This finishes the definition of $g_t$ for $t\in(0,\de)$. The reasoning
behind the definition -- why $\hat g$ is shrunk by a factor of $t^6$,
but the cut-off functions change between radii $t$ and $t^2$, for
instance -- will emerge in \S\ref{s24}, where we show that for this
definition the scalar curvature of $g_t$ is close to the constant
function $\nu$ in the $L^{n/2}$ norm.

\subsection{Combining two metrics of constant scalar curvature $\nu$}
\label{s22}

A family of metrics will now be defined on the connected sum of two
Riemannian manifolds with constant scalar curvature $\nu$. We shall do
this by writing down a zero scalar curvature Riemannian manifold with
two asymptotically flat ends, and gluing one end into each of the
constant scalar curvature manifolds using the method of \S\ref{s21};
the new manifold will form the `neck' in between.

Let $N$ be $\R^n\sm\{0\}$ with metric $g_{\sst N}=(1+\md{v}^{-(n-2)}
)^{p-2}h$. Then $N$ has two ends, and $g_{\sst N}$ is asymptotically
flat at each end. The involution $v\mapsto v/\ms{v}$ is an isometry.
Also, as $\De\md{v}^{-(n-2)}=0$, equation \eq{s1eq2} shows that
$g_{\sst N}$ has zero scalar curvature.

Let $(M',g')$ and $(M'',g'')$ be Riemannian $n$-manifolds with
constant scalar curvature $\nu$; applying homotheties if necessary
we shall assume that $\nu=1,0$ or $-1$. We shall use the gluing
method of \S\ref{s21} to glue the two asymptotically flat ends of
$(N,g_{\sst N})$ into $(M',g')$ and $(M'',g'')$. Let $M=M'\# M''$.
A family of metrics $\{g_t:t\in(0,\de)\}$ on $M$ will be defined,
such that $g_t$ resembles the union of $(M',g')$ and $(M'',g'')$
joined by a small `neck' of radius $t^6$, modelled upon $(N,t^{12}
g_{\sst N})$. It will be done briefly, as the treatment
generalizes~\S\ref{s21}.

Suppose that $M',M''$ contain points $m',m''$ with neighbourhoods in
which $g',g''$ are conformally flat. (In \S\ref{s34} this assumption
will be dropped.) Thus $M',M''$ contain open balls $B',B''$ with
diffeomorphisms $\Phi',\Phi''$ from $B_{\de}(0)$ in $\R^n$ to $B',B''$,
such that $\Phi'(0)=m',\Phi''(0)=m''$ and
$(\Phi')^*(g')=(\psi')^{p-2}h$, $(\Phi'')^*(g'')=(\psi'')^{p-2}h$, 
for some functions $\psi',\psi''$ on $B_{\de}(0)$. By choosing 
different conformal identifications with $B_{\de}(0)$ if necessary, 
we may suppose that $\psi'(0)=\psi''(0)=1$ and $\d\psi'(0)=\d\psi''(0)=0$, 
so that $\psi'(v)=1+O'(\ms{v})$ and~$\psi''(v)=1+O'(\ms{v})$.

For any $t\in(0,\de)$, define $M$ and the conformal class of $g_t$ by
\begin{equation*}
M=\Bigl(M'\sm\!\Phi'\bigl[\ovB_{t^2}(0)\bigr]\Bigr)\!\amalg\!
\Bigl(M''\sm\!\Phi''\bigl[\ovB_{t^2}(0)\bigr]\Bigr)\!\amalg\!
\bigl\{v\in N:t^5\!<\!\md{v}\!<\!t^{-5}\bigr\}\Bigl/\bigti_t,
\end{equation*} 
where $\bigti_t$ is the
equivalence relation defined by
\begin{gather*}
\Phi'[t^6v]\,\,\bigti_t\,\, v
\quad\text{if $v\in N$ and $t^{-4}<\md{v}<t^{-5}$, and}\\
\Phi''\left[\frac{t^6v}{\ms{v}}\right]\,\,\bigti_t\,\, v
\quad\text{if $v\in N$ and $t^{5}<\md{v}<t^{4}$.}
\end{gather*}

The conformal class $[g_t]$ of $g_t$ is then given by the restriction 
of the conformal classes of $g',g''$ and $g_{\sst N}$ to the open sets 
of $M',M''$ and $N$ that make up $M$; this definition makes sense 
because the conformal classes agree on the annuli of overlap where the
three open sets are glued by $\bigti_t$. Let $A_t$ be this region of 
gluing in $M$. Then $A_t$ is diffeomorphic via $\Phi'$ and $\Phi''$ 
to two copies of the annulus~$\{v\in\R^n:t^2<\md{v}<t\}$.

To define a metric $g_t$ within this conformal class, let $g_t=g'$ on
the component of $M\sm A_t$ coming from $M'$, $g_t=g''$ on the component 
of $M\sm A_t$ coming from $M''$, and $g_t=t^{12}g_{\sst N}$ on the 
component on $M\sm A_t$ coming from $N$. So it remains to choose a 
conformal factor on $A_t$ itself. Using $\Phi'$ and $\Phi''$, this
is the same as choosing a conformal factor for two copies of the subset
$\{v\in\R^n:t^2<\md{v}<t\}$ of~$\R^n$.

As in \S\ref{s21} define a partition of unity $\be_1,\be_2$ on $A_t$,
and define $\psi_t$ by $\psi_t(v)=\be_1(v)\psi'(v)+\be_2(v)(1+t^{6(n-2)}
\md{v}^{-(n-2)})$ on the component of $A_t$ coming from $M'$, and
$\psi_t(v)=\be_1(v)\psi''(v)+\be_2(v)(1+t^{6(n-2)}\md{v}^{-(n-2)})$
on the component coming from $M''$. Here $\psi_t$ is thought of as
a function on $A_t$, which is identified by $\Phi'$ and $\Phi''$
with two disjoint copies of $\{v\in\R^n:t^2<\md{v}<t\}$. Now let
$g_t$ be $\psi_t^{p-2}h$ in $A_t$, where $h$ is the push-forward
to $A_t$ by $\Phi',\Phi''$ of the standard metric on
$\{v\in\R^n:t^2<\md{v}<t\}$. This completes the definition of~$g_t$.

\subsection{The relationship between \S\ref{s21} and \S\ref{s22}
when $\nu=1$}
\label{s23}

Suppose that $(M',g')$ and $(M'',g'')$ are compact Riemannian
$n$-manifolds with scalar curvature 1, and $m',m''$ are points
in $M',M''$ with conformally flat neighbourhoods. Then \S\ref{s21}
and \S\ref{s22} with $\nu=1$ define two families of metrics
$\bigl\{g_t:t\in(0,\de)\bigr\}$ on $M=M'\# M''$. Reversing the
r\^oles of $M'$ and $M''$ in \S\ref{s21} we get a third family,
by shrinking the stereographic projection of $(M',g')$ from $m'$
by a factor $t^{12}$ and gluing it into $(M'',g'')$. What is the
relationship between these three families of metrics on~$M$?

Encoded in the choice of maps $\Phi',\Phi''$ and $\Xi''$ above is a
choice of isometric isomorphism $T_{m'}M'\cong T_{m''}M''$. Because
of the simplifying assumption of conformal flatness near $m',m''$ it
turns out that this isomorphism is the {\it only} arbitrary choice
involved in the construction of the conformal classes $[g_t]$, and
if we choose the same isomorphism $T_{m'}M'\cong T_{m''}M''$ in
the three cases, then the three families actually define the {\it same}
1-parameter family of conformal classes $\bigl\{[g_t]:t\in(0,\de)
\bigr\}$ for small~$t$.

However, the parametrization of the conformal classes by $t$ differs
in \S\ref{s21} and \S\ref{s22}. Identifying $M'$ and $M''$ with
$\R^n$ conformally near $m',m''$, and reasoning from the definition
of stereographic projection, we find that in \S\ref{s21} $\bigti_t$
identifies $v'\in M'$ with $v''\in M''$ only if $\md{v'}\cdot\md{v''}
=t^6$, but in \S\ref{s22} $\bigti_t$ identifies $v'\in M'$ with
$v''\in M''$ only if~$\md{v'}\cdot\md{v''}=t^{12}$.

Because of this, the conformal class $[g_t]$ constructed in \S\ref{s22}
is the same as the conformal class $[g_{t^2}]$ constructed in
\S\ref{s21}. We adopted this mildly inconvenient convention because with
it we will be able to prove results simultaneously for the metrics of
\S\ref{s21} and \S\ref{s22}, without changing the powers of $t$ involved.

\subsection{Estimating the scalar curvature of $g_t$}
\label{s24}

We now show that the scalar curvature of $g_t$ approaches the
constant value $\nu$ in the $L^{n/2}$ norm as~$t\ra 0$.

\begin{prop} Let\/ $\{g_t:t\in (0,\de)\}$ be one of the families
of metrics defined on $M=M'\# M''$ in \S\ref{s21} and\/ \S\ref{s22}.
Let the scalar curvature of\/ $g_t$ be $\nu-\ep_t$. Then there exist\/
$Y,Z>0$ such that\/ $\md{\ep_t}\le Y$ and\/ ${}^{g_t}\lnm{\ep_t}{n/2}\le
Zt^2$ for~$t\in(0,\de)$.
\label{p231}
\end{prop}

\begin{proof} The proof will be given for the metrics $g_t$ of
\S\ref{s21} only, the modifications for \S\ref{s22} being left to
the reader. We first derive an expression for $\ep_t$. Outside 
$A_t$, the scalar curvature of $g_t$ is $\nu$ and 0 on the regions
coming from $M'$ and $M''$. On $A_t$, calculating with \eq{s1eq2} gives
\e
\begin{split}
\nu-\ep_t(v)=\nu\be_1&(v)(\psi'(v))^{(n+2)/(n-2)}
\psi_t^{-(n+2)/(n-2)}(v)\\
&+a\psi_t^{-(n+2)/(n-2)}(v)\bigl(\De\be_1(v)\bigr)
\bigl(\psi'(v)-\xi(t^{-6}v)\bigr)\\
&-2a\psi_t^{-(n+2)/(n-2)}(v)\bigl(\na\be_1(v)\bigr)\cdot
\bigl(\na\bigl(\psi'(v)-\xi(t^{-6}v)\bigr)\bigr),
\end{split}
\label{s2eq2}
\e
since~$\be_1+\be_2=1$.

As $\psi'(v)=1+O'(\ms{v})$, $\xi(v)=1+O'(\md{v}^{2-n})$, and 
$t^2<\md{v}<t$, it follows easily that $\psi'(v)-\xi(t^{-6}v)=
O'(\ms{v})$. The reason for choosing to scale $\hat g$ by a factor of 
$t^{12}$ whilst making $\be_1$ change between $t^2$ and $t$ is to make 
this estimate work -- the first power has to be as high as 12 to work 
in dimension 3. Substituting it into \eq{s2eq2} gives that
\begin{align*}
\md{\ep_t}\le\md{\nu}\cdot
&\bmd{\be_1(\psi')^{(n+2)/(n-2)}\psi_t^{-(n+2)/(n-2)}-1}\\
+&\psi_t^{-(n+2)/(n-2)}\cdot\bigl\{
\md{\na\be_1}O(\md{v})+\md{\De\be_1}O(\ms{v})\bigr\}.
\end{align*}

Using a lower bound for $\psi'$ we find that on the subannulus
$t^2<\md{v}<t$, the estimate $\md{\psi_t(v)}\ge C_0>0$
holds for some $C_0$ and all $t\in(0,\de)$. Using this to get
rid of the $\psi_t$ terms on the r.h.s., and an upper bound for
$\psi'$, it can be seen that
\e
\md{\ep_t}=O(1)+\md{\na\be_1}O(\md{v})+\md{\De\be_1}
O\bigl(\ms{v}\bigr)
\label{s2eq3}
\e
on the subannulus $t^2<\md{v}<t$. But
\e
\md{\na\be_1}=\frac{\left\vert\frac{\d\tau_1}{\d x}\right\vert
}{\md{v\log t}}=O\bigl(\md{v}^{-1}\bigr),
\label{s2eq4}
\e
and in a similar way $\md{\De\be_1}=O\bigl(\md{v}^{-2}\bigr)$.
Substituting into \eq{s2eq3}, we find that for all $t\in(0,\de)$,
$\md{\ep_t}\le Y$ on the subannulus $t^2<\md{v}<t$, for
some~$Y\ge\md{\nu}$.

Thus $\md{\ep_t}\le Y$ on $A_t$, and outside $A_t$, $\ep_t=0$ on the
component coming from $M'$, and $\md{\ep_t}=\md{\nu}\le Y$ on the
component coming from $M''$. Therefore $\md{\ep_t}\le Y$, giving
the first part of the proposition. To prove the second part, observe
that by the estimates on $\psi_t$ above, the support of $\ep_t$ has
volume $\le C_1t^n$, for some $C_1>0$. So ${}^{g_t}\lnm{\ep_t}{n/2}\le
Zt^2$, where~$Z=YC_1^{2/n}$.
\end{proof}

\subsection{A uniform bound for a Sobolev embedding}
\label{s25}

If $(M,g)$ is a compact Riemannian $n$-manifold, then $L^2_1(M)$ is
continuously embedded in $L^{2n/(n-2)}(M)=L^p(M)$ by the {\it Sobolev
Embedding Theorem}, \cite[Th.~2.20]{Aubi}. This means that $L^2_1(M)
\subset L^p(M)$, and $\lnm{\phi}{p}\le A\snm{\phi}{2}{1}$ for all
$\phi\in L^2_1(M)$, and some $A>0$ depending on $g$. We shall prove
this holds for the metrics $g_t$ of \S\ref{s21} and \S\ref{s22},
with $A$ independent of~$t$.

\begin{prop} Let\/ $\{g_t:t\in(0,\de)\}$ be one of
the families of metrics defined on $M=M'\# M''$ in \S\ref{s21} and\/
\S\ref{s22}. Then there exist\/ $A>0$ and\/ $\ze\in(0,\de)$ such that\/
${}^{g_t}\lnm{\phi}{p}\le A\cdot{}^{g_t}\snm{\phi}{2}{1}$ whenever
$\phi\in L^2_1(M)$ and\/~$0<t\le\ze$.
\label{p241}
\end{prop}

\begin{proof} First consider the metrics of \S\ref{s21}. The proof
works by proving similar Sobolev embedding results on the component
manifolds $(M',g')$ and $(\hat M,t^{12}\hat g)$ that make up $(M,g_t)$,
and then `gluing' them together. These are given in the following
lemmas, the first following from the Sobolev Embedding
Theorem~\cite[Th.~2.20]{Aubi}.

\begin{lem} There exists $D_1>0$ with\/ ${}^{g'}\lnm{\phi}{p}\le
D_1\cdot{}^{g'}\snm{\phi}{2}{1}$ for all\/~$\phi\in L^2_1(M')$.
\label{l242}
\end{lem}

\begin{lem} There exists $D_2>0$ with\/ ${}^{t^{12}\hat g}
\lnm{\phi}{p}\le D_2\cdot{}^{t^{12}\hat g}\lnm{\na\phi}{2}$ for
all\/ $\phi\in L^2_1(\hat M)$ and\/~$t>0$.
\label{l243}
\end{lem}

\begin{proof} The inequality $\lnm{\phi}{p}\le D_2\lnm{\na\phi}{2}$
is invariant under homotheties, and so it is enough to prove the
lemma when $t=1$. Bartnik \cite[\S 1-\S 2]{Bart} reviews the theory
of weighted Sobolev spaces on asymptotically flat manifolds. We
shall apply \cite[Th.~1.2 \& Th.~1.3]{Bart}, which hold on
asymptotically flat manifolds by \cite[p.~676]{Bart}. Part (iv)
of \cite[Th.~1.2]{Bart} with $k=1$, $p=q=2$ and $\de=(2-n)/2$
shows that there exists $B>0$ such that 
\begin{equation*}
\left(\int_{\hat M}\md{\phi}^{2n/(n-2)}\d V_{\hat g}\right)^{(n-2)/n}\le
B\left(\int_{\hat M}\md{\na\phi}^2\d V_{\hat g}+\int_{\hat M}\md{\phi}^2
\si^{-2}\d V_{\hat g}\right)
\end{equation*}
for all $\phi\in L^2_1(N)$, where $\si:N\ra[1,\iy)$ is defined by
$\si(x)^2=1+d(x,x_0)^2$ for some base point $x_0\in N$. But part
(i) of \cite[Th.~1.3]{Bart} with $p=2$ and $\de=(2-n)/2$ shows
that there exists $C>0$ such that
\begin{equation*}
\int_{\hat M}\md{\phi}^2\si^{-2}\d V_{\hat g}\le
C\int_{\hat M}\md{\na\phi}^2\d V_{\hat g}.
\end{equation*}
Combining the last two equations and taking square roots proves
the lemma, with~$D_2=(B(1+C))^{1/2}$.
\end{proof}

The volume form of $g_t$ on $A_t$ is $\psi_t^p\,\d V_h$, so the
contribution from $A_t$ to $\int_M\phi^p\,\d V_{g_t}$ is $\int_{A_t}
(\psi_t\phi)^p\,\d V_h$. Using \eq{s2eq1} we may eliminate $\psi_t$,
divide into integrals on the component manifolds, and show that
${}^{g_t}\lnm{\phi}{p}\le{}^{g'}\lnm{\be_1\phi}{p}+{}^{t^{12}\hat g}
\lnm{\be_2\phi}{p}$ for $\phi\in L^2_1(M)$. Now
\begin{align*}
{}^{g'}\snm{\be_1\phi}{2}{1}&\le
{}^{g'}\lnm{\be_1\phi}{2}+{}^{g'}\lnm{\be_1\na\phi}{2}+
{}^{g'}\blnm{\phi\md{\na\be_1}}{2}\\
&\le {}^{g'}\lnm{\be_1\phi}{2}+{}^{g'}\lnm{\be_1\na\phi}{2}+
{}^{g'}\blnm{\phi\vert_{A_t}}{p}
\cdot{}^{g'}\blnm{\d\be_1\vert_{A_t}}{n},
\end{align*}
by H\"older's inequality. But as $\lnm{\d\be_1\vert_{A_t}}{n}$
is a conformally invariant norm we have ${}^{g'}\blnm{\d\be_1
\vert_{A_t}}{n}={}^{h}\blnm{\d\be_1\vert_{A_t}}{n}$, where
$h$ is the standard metric on~$\R^n$.

Combining the last few equations and the obvious analogues for
$\hat M$, and remembering that $\d\be_1+\d\be_2=0$, we obtain
\begin{align*}
{}^{g_t}\lnm{\phi}{p}&\le
D_1\cdot{}^{g'}\snm{\be_1\phi}{2}{1}+D_2\cdot{}^{t^{12}\hat g}
\snm{\be_2\phi}{2}{1}\\
&\le D_1\bigl({}^{g'}\lnm{\be_1\phi}{2}+
{}^{g'}\lnm{\be_1\na\phi}{2}+
{}^{g'}\blnm{\phi\vert_{A_t}}{p}
\cdot{}^{h}\blnm{\d\be_1\vert_{A_t}}{n}\bigr)\\
&+D_2\bigl({}^{t^{12}\hat g}\lnm{\be_2\phi}{2}+
{}^{t^{12}\hat g}\lnm{\be_2\na\phi}{2}+
{}^{t^{12}\hat g}\blnm{\phi\vert_{A_t}}{p}
\cdot{}^{h}\blnm{\d\be_1\vert_{A_t}}{n}\bigr)\\
&\le D_3\cdot{}^{g_t}\snm{\phi}{2}{1}
+D_4\cdot{}^{g_t}\lnm{\phi}{p}\cdot{}^{h}\blnm{\d\be_1\vert_{A_t}}{n},
\end{align*}
for all $t\in(0,\de)$ and some $D_3,D_4>0$ depending on $D_1,D_2$
and bounds for the ratios of the various conformal factors on~$A_t$.

From the last equation we see that Proposition \ref{p241} holds with
$A=2D_3$ provided $D_4\cdot{}^{h}\blnm{\d\be_1\vert_{A_t}}{n}\le\ha$.
As $\be_1$ depends only on $r=\md{v}$ we have
\begin{equation*}
{}^{h}\lnm{\d\be_1}{n}^n=\om_{n-1}\cdot\int_{t^2}^t
r^{n-1}\Bigl\vert\frac{\d\be_1}{\d r}\Bigr\vert^n\d r
=\om_{n-1}\md{\log t}^{1-n}\cdot\int_1^2
\Bigl\vert\frac{\d\tau_i}{\d x}\Bigr\vert^n\d x,
\end{equation*}
where $\om_{n-1}$ is the volume of the unit sphere ${\cal S}^{n-1}$
in $\R^n$. So when $t$ is small, ${}^h\lnm{\d\be_1}{n}$ is small,
and there exists $\ze\in(0,\de)$ such that $D_4{}^{h}\blnm{\d\be_1
\vert_{A_t}}{n}\le\ha$ when $0<t\le\ze$.
This completes the proposition for the metrics of \S\ref{s21}.
The proof for the metrics of \S\ref{s22} requires only simple
modifications, and we leave it to the reader.
\end{proof}

\section{Existence results for scalar curvature $\pm 1$}
\label{s3}

Let $M$ be the manifold of \S\ref{s21} or \S\ref{s22} with one of the metrics $g_t$
defined there, and denote its scalar curvature by $S$. As in \S\ref{s12}, a
conformal change to $\ti g_t=\psi^{p-2}g_t$ may be made, and the condition
for $\ti g_t$ to have constant scalar curvature $\nu$ is the Yamabe
equation
\e
a\De\psi+S\psi=\nu\md{\psi}^{p-1}.
\label{s3eq1}
\e
Now the metrics $g_t$ have scalar curvature close to $\nu$, so let
$S=\nu-\ep$; then by Proposition \ref{p231}, ${}^{g_t}\lnm{\ep}{n/2}\le Zt^2$. Also
we would like the conformal change to be close to 1, so put $\psi=1+\phi$,
where we aim to make $\phi$ small. Substituting both of these changes into
\eq{s3eq1} gives
\e
a\De\phi-\nu b\phi=\ep+\ep\phi+\nu f(\phi),
\label{s3eq2}
\e
where $b=4/(n-2)$ and $f(t)=\md{1+t}^{(n+2)/(n-2)}-1-(n+2)t/(n-2)$.

In this section, we shall suppose that $\nu=\pm 1$, as the zero scalar
curvature case requires different analytic treatment and will be
considered in \S\ref{s4}. Equation \eq{s3eq2} has been written so
that on the left is a linear operator $a\De-\nu b$ applied to $\phi$,
and on the right are the `error terms'.

The method of \S\ref{s31} is to define by induction a sequence
$\{\phi_i\}_{i=0}^\iy$ of functions in $L^2_1(M)$ by $\phi_0=0$, and
\begin{equation*}
a\De\phi_i-\nu b\phi_i=\ep+\ep\phi_{i-1}+\nu f(\phi_{i-1}).
\end{equation*}
This depends upon being able to invert the operator $a\De-\nu b$,
and we consider the existence and size of the inverse in \S\ref{s32}
and \S\ref{s33}. Given this invertibility, we show that if $\ep$ is
sufficiently small, $\{\phi_i\}_{i=0}^\iy$ converges to $\phi\in
L^2_1(M)$ which is a weak solution of \eq{s3eq2}. Finally we show
that $\phi$ is smooth and $\psi=1+\phi$ is positive, so that $\ti g_t$
has constant scalar curvature~$\nu$.

In \S\ref{s32} and \S\ref{s33} we state the existence theorems
for constant positive and negative scalar curvature respectively
on connected sums, the main results of this section. Note that
\S\ref{s33} produces {\it three distinct} metrics of scalar curvature
1 in the conformal class of each suitable connected sum of manifolds
with scalar curvature 1, in contrast to the negative scalar curvature
case, where any metric of scalar curvature $-1$ is unique in its
conformal class.

\subsection{An existence result for constant scalar curvature}
\label{s31}

Fix $\nu=\pm 1$, and suppose $(M,g)$ is a compact Riemannian $n$-manifold.
Let $A,B,X$ and $Y$ be positive constants, to be chosen later. We write
down four properties, which $(M,g)$ may or may not satisfy:
\medskip

\noindent{\bf Property 1.} {\it The volume of\/ $M$ satisfies
$X/2\le\vol(M)\le X$.}
\medskip

\noindent{\bf Property 2.} {\it Let the scalar curvature of\/ $g$
be $\nu-\ep$. Then $\md{\ep}\le Y$.}
\medskip

\noindent{\bf Property 3.} {\it Whenever $\phi\in L^2_1(M)$,\ $\phi\in
L^p(M)$, and\/ $\lnm{\phi}{p}\le A\snm{\phi}{2}{1}$.}
\medskip

\noindent{\bf Property 4.} {\it For every $\xi\in L^{2n/(n+2)}(M)$, there
exists a unique $\phi\in L^2_1(M)$ such that\/ $a\De\phi-\nu b\phi=\xi$
holds weakly. Moreover, $\snm{\phi}{2}{1}\le B\lnm{\xi}{2n/(n+2)}$.}
\medskip

We of course think of $(M,g)$ as being the manifold $M$ of \S\ref{s21} or
\S\ref{s22}, with one of the metrics $g_t$ defined there. Then Property 1 is
clear from the definitions, Property 2 comes from Proposition \ref{p231},
Property 3 from Proposition \ref{p241}, and Property 4 remains to be
proved. In terms of these properties, we state the next result, which
is the core of the analysis of this section.

\begin{thm} Let\/ $A,B,X,Y>0$ and\/ $n\ge 3$ be given. Then there exist\/
$W,c>0$ depending only upon $A,B,X,Y$ and $n$, such that if\/ $(M,g)$
satisfies Properties $1$--$4$ above and\/ $\lnm{\ep}{n/2}\le c$, then
$g$ admits a smooth conformal rescaling to $\ti g=(1+\phi)^{p-2}g$,
with constant scalar curvature $\nu$, and\/~$\snm{\phi}{2}{1}\le
W\lnm{\ep}{n/2}$.
\label{t311}
\end{thm}

\begin{proof} Suppose that $(M,g)$ satisfies Properties 1--4 above.
Define a map $T:L^2_1(M)\ra L^2_1(M)$ by $T\eta=\xi$, where
\e
a\De\xi-\nu b\xi=\ep+\ep\eta+\nu f(\eta).
\label{s3eq3}
\e

By Property 4, $\xi$ exists and is unique, provided that the right hand
side is in $L^{2n/(n+2)}(M)$. So it must be shown that if $\eta\in L^2_1(M)$,
then $\ep+\ep\eta+\nu f(\eta)\in L^{2n/(n+2)}(M)$. Now $\ep\in
L^{n/2}(M)$ implies that $\ep\in L^{2n/(n+2)}(M)$; by the Sobolev
embedding theorem $\eta\in L^{2n/(n-2)}(M)$, and as $\ep\in L^{n/2}(M)$
it follows that $\ep\eta\in L^{2n/(n+2)}(M)$. Thus the first two terms
are in $L^{2n/(n+2)}(M)$. For the third term, as $\eta\in L^{2n/(n-2)}(M)$,
$1+\eta$ is too, and $(1+\eta)^{(n+2)/(n-2)}\in L^{2n/(n+2)}(M)$. This deals
with the first part of $f(\eta)$, and the last two parts are trivially in
$L^{2n/(n+2)}(M)$. Therefore the right hand of \eq{s3eq3} is in
$L^{2n/(n+2)}(M)$, and the map $T$ is well defined.

Now define a sequence $\{\phi_i\}_{i=0}^\iy$ of elements of $L^2_1(M)$ by
$\phi_i=T^i(0)$. Our first goal is to prove that if $\lnm{\ep}{n/2}$ is
sufficiently small, then this sequence converges in $L^2_1(M)$. Setting
$\phi$ as the limit of the sequence, \eq{s3eq3} implies that $\phi$ will
satisfy \eq{s3eq2}, as we would like. This will be achieved via the next
lemma.

\begin{lem} Suppose $L$ is a Banach space with norm $\nm{.}$, 
and\/ $T:L\ra L$ is a map satisfying $\nm{T(0)}\le F_0s$ and
\e
\begin{split}
\nm{T(v)-T(w)}\le\,&\nm{v-w}\Bigl\{F_1s+F_2\bigl(\nm{v}
+\nm{w}\bigr)\\
&+F_3\bigl(\nm{v}^{4/(n-2)}+\nm{w}^{4/(n-2)}
\bigr)\Bigr\}
\end{split}
\label{s3eq4}
\e
for all\/ $v,w\in L$, where $F_0,F_1,F_2,F_3,s>0$. Then there exists
$W>0$ depending only on $F_0,F_1,F_2,F_3$ and\/ $n$, such that if\/
$s$ is sufficiently small, then the sequence $\{\phi_i\}_{i=0}^\iy$
defined by $\phi_i=T^i(0)$ converges to a limit\/ $\phi$ in $L$,
satisfying~$\nm{\phi}\le Ws$.
\label{l312}
\end{lem}

\begin{proof} Putting $v=\phi_{i-1}$ and $w=0$ into \eq{s3eq4} gives 
$\nm{\phi_i-T(0)}\le F_1s\nm{\phi_{i-1}}+F_2\nm{\phi_{i-1}}^2
+F_3\nm{\phi_{i-1}}^{(n+2)/(n-2)}$, and as $\nm{T(0)}\le F_0s$ this
implies that $\nm{\phi_i}\le \chi(\nm{\phi_{i-1}})$, where
$\chi(x)=F_0s+F_1sx+F_2x^2+F_3x^{(n+2)/(n-2)}$.

From the form of this equation, it is clear that there exists $W>0$
depending only on $F_0,F_1,F_2,F_3,n$ such that when $s$ is small,
there exists an $x$ with $0<x\le Ws$ and $2\chi(x)=x$. Suppose $s$
is small enough, so that such an $x$ exists. Now $\phi_0=0$, so that
$\nm{\phi_0}\le x$, and if $\nm{\phi_{i-1}}\le x$ then $\nm{\phi_i}\le
\chi(\nm{\phi_{i-1}})\le\chi(x)\le x$. Thus by induction,
$\nm{\phi_i}\le x$ for all~$i$.

Put $v=\phi_i$ and $w=\phi_{i-1}$ in \eq{s3eq4}. This gives
$\nm{\phi_{i+1}-\phi_i}\le\nm{\phi_i-\phi_{i-1}}\cdot
\bigl(F_1s+2F_2x+2F_3x^{4/(n-2)}\bigr)$, using the inequality $\nm{\phi_i}\le
x$ that we have just proved. Dividing the equation $x=2\chi(x)$ by $x$ and
subtracting some terms it follows that $1>F_1s+2F_2x+2F_3x^{4/(n-2)}>0$, and
so $\{\phi_i\}_{i=0}^\iy$ converges, by comparison with a geometric
series. Let the limit of the sequence be $\phi$. Then as $\nm{\phi_i}\le x\le
Ws$ for all~$i$, by continuity $\phi$ also satisfies $\nm{\phi}\le Ws$.
\end{proof}

To apply the lemma we must show that $T:L^2_1(M)\ra L^2_1(M)$ defined
above satisfies the hypotheses. Let $s=\lnm{\ep}{n/2}$. We will define
$F_0,F_1,F_2,F_3>0$ depending only on $A,B,X$ and $n$, such that 
\eq{s3eq4} holds.

Putting $\eta=0$ in \eq{s3eq3}, and applying Properties 1 and 4, 
we see that 
\begin{equation*}
\snm{T(0)}{2}{1}\le B\lnm{\ep}{2n/(n+2)}\le B\vol(M)^{(n-2)/2n}
\lnm{\ep}{n/2}\le BX^{(n-2)/2n}s, 
\end{equation*}
so let $F_0=BX^{(n-2)/2n}$. From the
definition of $f$, it can be seen that
\e
\bmd{f(x)-f(y)}\le\bmd{x-y}\cdot
\Bigl(F_4\bigl(\md{x}+\md{y}\bigr)+F_5\bigl(\md{x}^{4/(n-2)}+
\md{y}^{4/(n-2)}\bigr)\Bigr),
\label{s3eq5}
\e
where $F_4,F_5$ are constants depending only on $n$, and $F_4=0$ if
$n\ge 6$. Let $\eta_1,\eta_2\in L^2_1(M)$, and let $T(\eta_i)=\xi_i$. Then
taking the difference of \eq{s3eq3} with itself for $i=1,2$ we get
\begin{equation*}
(a\De-\nu b)(\xi_1-\xi_2)=\ep\cdot(\eta_1-\eta_2)+
\nu\bigl(f(\eta_1)-f(\eta_2)\bigl).
\end{equation*}

Applying Property 4 and making various estimates gives
\begin{align*}
\snm{\xi_1-\xi_2}{2}{1}&\le 
B\Bigl(\lnm{\ep\cdot(\eta_1-\eta_2)}{2n/(n+2)}+
\md{\nu}\blnm{f(\eta_1)-f(\eta_2)}{2n/(n+2)}\Bigr)\\
&\le\snm{\eta_1-\eta_2}{2}{1}\cdot
\Bigl(F_1s+F_2\bigl(\snm{\eta_1}{2}{1}+\snm{\eta_2}{2}{1}\bigr)\\
&\qquad\qquad\qquad\quad +F_3\bigl(\snm{\eta_1}{2}{1}^{4/(n-2)}
+\snm{\eta_2}{2}{1}^{4/(n-2)}\bigr)\Bigr),
\end{align*}
where $F_1=AB$,\ $F_2=A^2BF_4X^{(6-n)/2n}$ and
$F_3=A^{(n+2)/(n-2)}BF_5$. The calculation uses H\"older's inequality,
\eq{s3eq5}, Properties 1 and 3, the expression $\lnm{\eta}{r}\le
\lnm{\eta}{s}(\vol M)^{(s-r)/rs}$ when $1\le r<s$ and
$\eta\in L^r(M) \subset L^s(M)$, and the fact that $F_4=0$ if~$n\ge 6$.

This inequality is \eq{s3eq4} for the operator $T:L^2_1(M)\ra L^2_1(M)$.
So putting $L=L^2_1(M)$ and applying Lemma \ref{l312}, there is a constant
$W>0$ depending only on $F_0,F_1,F_2,F_3$ and $n$, such that if
$\lnm{\ep}{n/2}$ is sufficiently small, then the sequence
$\{\phi_i\}_{i=0}^\iy$ defined by $\phi_i=T^i(0)$ converges to a
limit $\phi$, satisfying $\nm{\phi}\le W\lnm{\ep}{n/2}$. Now $W$
depends only on $F_0,\dots,F_3$ and $n$, and these depend only on $n,A,B$ and
$X$, so $W$ depends only on $n,A,B$ and $X$. Since $\phi_i=T(\phi_{i-1})$ and
$T$ is continuous, taking the limit gives $\phi=T(\phi)$, so \eq{s3eq3} shows
that $\phi$ satisfies \eq{s3eq2} weakly. Thus we have proved the following lemma:

\begin{lem} There exists $W>0$ depending only on $n,A,B$ and\/ $X$,
such that if\/ $\lnm{\ep}{n/2}$ is sufficiently small, then there
exists $\phi\in L^2_1(M)$ satisfying \eq{s3eq2} weakly, with\/
$\snm{\phi}{2}{1}\le W\lnm{\ep}{n/2}$.
\label{l313}
\end{lem}

Thus weak solutions $\phi$ of \eq{s3eq2} do exist for small
$\lnm{\ep}{n/2}$, and for these $\psi=1+\phi$ is a weak solution of
\eq{s1eq3}. But for $\ti g=\psi^{p-2}g$ to be a metric we need $\psi$
to be smooth and positive. Proposition \ref{p122} shows that $\psi\in
C^2(M)$, and is $C^\iy$ wherever it is nonzero, so it remains to show
that~$\psi>0$.

Examples of manifolds (with negative scalar curvature) can be found for
which \eq{s1eq3} admits solutions that change sign, so there is something
to be proved. This difficulty does not arise in the proof of the Yamabe
problem, as there $\psi$ is the limit of a minimizing sequence of
positive functions, so $\psi\ge 0$ automatically.

We deal with this problem in the following proposition, by proving that if
$\psi=1+\phi$ is a solution to \eq{s1eq3} that is negative somewhere, then
$\snm{\phi}{2}{1}$ must be at least a certain size. So if $\phi$ is small
in $L^2_1(M)$, then $\psi=1+\phi\ge 0$. We then show that $\psi>0$ using a
maximum principle.

\begin{prop} If\/ $\snm{\phi}{2}{1}$ is sufficiently small, then~$\psi\ge 0$.
\label{p414}
\end{prop}

\begin{proof} Let $\xi=\min(\psi,0)$. Then $\xi\in L^2_1(M)$ and 
$\int_M\xi\De\psi\,\d V_g=\int_M\ms{\na\xi}\,\d V_g$, since $\xi\De\psi
=\ms{\na\xi}+\ha\De\xi^2$, and $\int_M\De\xi^2\,\d V_g=0$. So multiplying
\eq{s1eq3} by $\xi$ and integrating over $M$ gives
\e
\int_M(a\ms{\na\xi}+S\xi^2+\nu\md{\xi}^{p})\,\d V_g=0,
\label{s3eq6}
\e
as $\xi=-\md{\xi}$. Also, from Property 3 and H\"older's inequality,
we have
\e
A\snm{\xi}{2}{1}\ge\lnm{\xi}{p}\ge\vol(\supp\xi)^{-1/n}\lnm{\xi}{2}.
\label{s3eq7}
\e

Using \eq{s3eq7} to eliminate $\int_Ma\ms{\na\xi}\d V_g$ in \eq{s3eq6},
and using $S\ge\nu-Y$, $\md{\nu}=1$, we find that
\e
\lnm{\xi}{p}^p\ge(F_6\vol(\supp\xi)^{-2/n}-F_7)\lnm{\xi}{2}^2,
\label{s3eq8}
\e
for $F_6,F_7>0$ depending on $A,Y$ and $n$. Now $\lnm{\xi}{2}\ge
\vol(\supp\xi)^{1/n}\lnm{\xi}{p}$ as above, and so either
$\lnm{\xi}{p}=0$, or we may substitute this in to \eq{s3eq8} to get
\e
\lnm{\xi}{p}^{4/(n-2)}\ge F_6-F_7\vol(\supp\xi)^{2/n}.
\label{s3eq9}
\e

If this holds then either $\lnm{\xi}{p}^{4/(n-2)}\ge F_6/2$ or 
$\vol(\supp\xi)\ge (F_6/2F_7)^{n/2}$. Both imply that $\lnm{\phi}{p}$
is bounded below by a positive constant, and by Property 3,
$\snm{\phi}{2}{1}$ is too. Conversely, if $\snm{\phi}{2}{1}$ is
smaller than this constant, then \eq{s3eq9} cannot hold and so
$\lnm{\xi}{p}=0$, which implies~$\psi\ge 0$.
\end{proof}

We are now ready to define the constant $c$ in Theorem \ref{t311}.
Let $c$ be small enough that three conditions hold: firstly,
$\lnm{\ep}{n/2}\le c$ implies $\lnm{\ep}{n/2}$ is sufficiently small
to satisfy Lemma \ref{l313}, so thata $\phi$ exists and satisfies
$\snm{\phi}{2}{1}\le cW$; secondly, that $\snm{\phi}{2}{1}\le cW$
implies $\snm{\phi}{2}{1}$ is sufficiently small to satisfy
Proposition \ref{p414}, so that $\psi=1+\phi\ge 0$; and thirdly,
that $\snm{\phi}{2}{1}\le cW$ implies $\snm{\phi}{2}{1}$ is
sufficiently small that $\phi$ cannot be the constant $-1$.
(As $X/2\le\vol(M)$ by Property 1, this depends only on~$X$.)

Then $c$ depends only on $n,A,B,X$ and $Y$, as the three conditions each
separately do. Thus if $\lnm{\ep}{n/2}\le c$, then there exists $\phi$
with $\snm{\phi}{2}{1}\le W\lnm{\ep}{n/2}$, such that $\psi=1+\phi\ge 0$
and satisfies \eq{s1eq3}. By the third condition on $c$, $\psi$ is not
identically zero. By Proposition \ref{p122}, $\psi\in C^2(M)$, and is
$C^\iy$ wherever it is nonzero. It remains to show that $\psi>0$
for $\ti g=\psi^{p-2}g$ to be nonsingular and have constant scalar
curvature $\nu$. This we achieve using the {\it strong maximum
principle}~\cite[Th.~2.6]{LePa}:

\begin{thm} Suppose $h$ is a nonnegative, smooth function on a 
connected manifold\/ $M$, and\/ $u\in C^2(M)$ satisfies $(\De+h)u\ge 0$. 
If\/ $u$ attains its minimum $m\le 0$, then $u$ is constant on~$M$.
\label{t315}
\end{thm}

As $M$ is compact and $S$ and $\psi$ are continuous, they are bounded on $M$,
and there is a constant $h\ge 0$ such that $S-\nu\psi^{p-2}\le h$ on $M$. Now
$M$ is connected, and $\psi\in C^2(M)$ satisfies \eq{s1eq3} and is nonnegative,
so $\psi$ satisfies $a\De\psi+h\psi\ge 0$. Thus by the strong maximum
principle, if $\psi$ attains the minimum value zero, then $\psi$ is
identically zero on $M$. But it has already been shown that this is not the
case, so $\psi$ cannot be zero anywhere and must be strictly positive. The
proof of Theorem \ref{t311} is therefore complete. 
\end{proof}

\subsection{Constant negative scalar curvature}
\label{s32}

We now construct metrics of scalar curvature $-1$ on connected sums
using the results of \S\ref{s31}. Fix $\nu=-1$, and consider the
metrics $g_t$ of \S\ref{s21} and \S\ref{s22}. Properties 1--3 of
\S\ref{s31} have already been dealt with, so it remains to show
that Property 4 holds for the metrics $g_t$, and that
${}^{g_t}\lnm{\ep_t}{n/2}$ is small when $t$ is small. As $\nu=-1$
Property 4 is about the invertibility of $a\De+b$, which is simple
as the eigenvalues of $\De$ are nonnegative.

\begin{lem} Let\/ $\{g_t:t\in (0,\de)\}$ be one of the families of 
metrics defined on $M=M'\# M''$ in \S\ref{s21} or \S\ref{s22}, and
let\/ $A,\ze$ be as in Proposition \ref{p241}. Then for all\/
$t\in(0,\ze]$ and\/ $\xi\in L^{2n/(n+2)}(M)$, there exists a unique
$\phi\in L^2_1(M)$ such that\/ $a\De\phi+b\phi=\xi$ holds weakly,
and\/ ${}^{g_t}\snm{\phi}{2}{1}\le B\,{}^{g_t}\lnm{\xi}{2n/(n+2)}$,
where~$B=A/b$.
\label{l321}
\end{lem}

\begin{proof} As $\De$ is self-adjoint and all its eigenvalues are
nonnegative, and as $a,b>0$, by some well-known analysis $a\De+b$ has a
right inverse, $T$ say, from $L^2(M)\ra L^2(M)$. Now $M$ is compact,
and so $L^2(M)\subset L^{2n/(n+2)}(M)$. Let $\xi\in L^2(M)$. We may
define $\phi\in L^2(M)$ by $\phi=T\xi$, and $a\De\phi+b\phi=\xi$ will
hold weakly.

It must first be shown that $\phi\in L^2_1(M)$ and that it satisfies the
inequality. Since $\phi,\xi\in L^2(M)$, $\int_M\phi\xi\d V_{g_t}$ exists,
and by subtraction $\int_M\phi\De\phi\d V_{g_t}$ exists as well. This is
$\int_M\ms{\na\phi}\d V_{g_t}$, and so $\phi\in L^2_1(M)$ by definition.

Multiplying the expression above by $\phi$ and integrating gives
$a\int_M\ms{\na\phi}\d V_{g_t}+b\,{}^{g_t}\lnm{\phi}{2}^2=\int_M\phi\xi\d
V_{g_t}$. As $a>b$, the l.h.s.\ is at least $b\,{}^{g_t}\snm{\phi}{2}{1}^2$,
and the r.h.s.\ is at most ${}^{g_t}\lnm{\phi}{p}{}^{g_t}
\lnm{\xi}{2n/(n+2)}$ by H\"older's inequality, since
$\phi\in L^p(M)$ by the Sobolev embedding theorem.
But ${}^{g_t}\lnm{\phi}{p}\le A\,{}^{g_t}\snm{\phi}{2}{1}$ by
Proposition \ref{p241}. Putting all this together gives
$b\,{}^{g_t}\snm{\phi}{2}{1}^2\le A\,{}^{g_t}\lnm{\xi}{2n/(n+2)}
{}^{g_t}\snm{\phi}{2}{1}$, and dividing by $b\,{}^{g_t}\snm{\phi}{2}{1}$
gives~${}^{g_t}\snm{\phi}{2}{1}\le B\,{}^{g_t}\lnm{\xi}{2n/(n+2)}$.

So far we have worked with $\xi\in L^2(M)$ rather than $L^{2n/(n+2)}(M)$. It
has been shown that the operator $T:L^2(M)\subset L^{2n/(n+2)}(M)\ra
L^2_1(M)$ is linear and continuous with respect to the $L^{2n/(n+2)}$ norm on
$L^2(M)$, and bounded by $B$. But therefore, by elementary functional
analysis, the operator $T$ extends uniquely to a continuous operator on the
closure of $L^2(M)$ in $L^{2n/(n+2)}(M)$, that is, $L^{2n/(n+2)}(M)$ itself.
Call this extended operator $\ov T$. Then for $\xi\in L^{2n/(n+2)}(M)$,
$\phi=\ov T\xi$ is a well-defined element of $L^2_1(M)$, satisfies
${}^{g_t}\snm{\phi}{2}{1}\le B\,{}^{g_t}\lnm{\xi}{2n/(n+2)}$, and $a\De\phi+b\phi=\xi$ holds
in the weak sense, by continuity. This concludes the proof. 
\end{proof}

All the previous work now comes together to prove the following two 
existence theorems for metrics of scalar curvature $-1$:

\begin{thm} Let\/ $(M',g')$ and\/ $(M'',g'')$ be compact Riemannian
$n$-manifolds with scalar curvature $-1$ and\/ $1$ respectively.
Suppose $M',M''$ contain points $m',m''$ with neighbourhoods in
which $g',g''$ are conformally flat.

As in \S\ref{s21}, define the family $\{g_t:t\in (0,\de)\}$ of
metrics on $M=M'\# M''$. Then there exists $C>0$ such that\/ $g_t$
admits a smooth conformal rescaling to $\ti g_t=(1+\phi)^{p-2}g_t$
with scalar curvature $-1$ for small\/ $t$,
and\/~${}^{g_t}\snm{\phi}{2}{1}\le Ct^2$.
\label{t322}
\end{thm}

\begin{thm} Let\/ $(M',g')$ and\/ $(M'',g'')$ be compact 
Riemannian $n$-manifolds with scalar curvature $-1$. Suppose
$M',M''$ contain points $m',m''$ with neighbourhoods in which
$g',g''$ are conformally flat.

As in \S\ref{s22}, define the family $\{g_t:t\in (0,\de)\}$ of metrics
on $M=M'\# M''$. Then there exists $C>0$ such that\/ $g_t$ admits a
smooth conformal rescaling to $\ti g_t=(1+\phi)^{p-2}g_t$ with scalar
curvature $-1$ for small\/ $t$, and\/~${}^{g_t}\snm{\phi}{2}{1}\le Ct^2$.
\label{t323}
\end{thm}

The proofs of the theorems are nearly the same, so only the first 
will be given. To get the second proof, change $\vol(M')$ to
$\vol(M')+\vol(M'')$ in the definition of~$X$.

\begin{proof}[Theorem \ref{t322}] Applying Propositions \ref{p231} and
\ref{p241} to the family $\{g_t:t\in (0,\de)\}$ gives a constant $Y$ for
Property 2 of \S\ref{s31}, and constants $A,\ze$ such that if $t\le\ze$ then
Property 3 holds for $g_t$ with constant $A$. By Lemma \ref{l321}, there is a
constant $B$ such that Property 4 also holds for $g_t$ when $t\le\ze$. 

It is clear that as $t\ra 0$, $\vol(M,g_t)\ra\vol(M')>0$. So there is 
a constant $X>0$ such that $X/2\le\vol(M,g_t)\le X$ for small enough $t$. 
This gives Property 1. Thus there are constants $n,A,B,X,Y$ such that 
Properties 1--4 of \S\ref{s31} hold for $(M,g_t)$ when $t$ is small.
Theorem \ref{t311} therefore gives a constant $c$ such that if
${}^{g_t}\lnm{\ep_t}{n/2}\le c$, we have the smooth conformal rescaling to 
a constant scalar curvature metric that we want.

But by Proposition \ref{p231}, ${}^{g_t}\lnm{\ep_t}{n/2}\le Zt^2$. So for small
enough $t$,\ ${}^{g_t}\lnm{\ep_t}{n/2}\le c$, and there exists a smooth conformal
rescaling to a metric $\ti g_t=(1+\phi)^{p-2}g_t$ which has scalar
curvature $-1$. Moreover, $\snm{\phi}{2}{1}\le W{}^{g_t}\lnm{\ep_t}{n/2}\le
WZt^2$, where $W$ is the constant given by Theorem \ref{t311}. Therefore 
putting $C=WZ$ completes the proof.
\end{proof}

\subsection{Constant positive scalar curvature}
\label{s33}

Now we construct metrics of scalar curvature 1 on connected sums. The
problems we encounter are in proving Property 4 of \S\ref{s31}, which now deals
with the invertibility of $a\De-b$, and they arise because $a\De$ may
have eigenvalues close to $b$. Our strategy is to show that if $a\De$ has
no eigenvalues in a fixed neighbourhood of $b$ on the component manifolds of
the connected sum, then for small $t$, $a\De$ has no eigenvalues in a
smaller neighbourhood of $b$ on~$(M,g_t)$. 

This is the content of the next theorem. We shall indicate here why the
theorem holds, but the proof we leave until the appendix, because it forms
a rather long and involved diversion from the main thread of the paper.

\begin{thm} Let\/ $\{g_t:t\in(0,\de)\}$ be one of the families of
metrics defined on $M=M'\# M''$ in \S\ref{s21} or \S\ref{s22},
and suppose that for some $\ga>0$, $a\De$ has no eigenvalues in 
$(b-2\ga,b+2\ga)$ on $(M',g')$ in \S\ref{s21}, and on both\/ $(M',g')$ 
and\/ $(M'',g'')$ in \S\ref{s22}. Then $a\De$ has no eigenvalues in
$(b-\ga,b+\ga)$ on $(M,g_t)$ for small\/~$t$.
\label{t331}
\end{thm}

Here is a sketch of the proof. Suppose $\phi$ is an eigenvector of
$a\De$ on $(M,g_t)$ for small $t$. Restricting $\phi$ to the portions
of $M$ coming from $M'$ and $M''$ and smoothing off gives functions on
$M',M''$. We try to show that one of these is close to an eigenvector
of $a\De$ on $M'$ or $M''$. This can be done except when $\phi$ is
large on the neck compared to the rest of the manifold.

But as the neck is a small region when $t$ is small, for $\phi$ to
be large there and small elsewhere means that $\phi$ must change
quickly around the neck, so that $\int_M\ms{\na\phi}\d V_{g_t}$ has
to be large compared to $\int_M\phi^2\d V_{g_t}$. Thus the eigenvalue
of $\phi$ must be large. Conversely, if the eigenvalue of $\phi$ is
close to $b$, then $\phi$ cannot be large on the neck compared to the
rest of $M$, and therefore either $M'$ or $M''$ must also have an
eigenvalue close to~$b$.

Using this result, Property 4 of \S\ref{s31} can be proved for the metrics:

\begin{lem} Let\/ $\{g_t:t\in(0,\de)\}$ be one of the families 
of metrics defined on $M=M'\# M''$ in \S\ref{s21} or \S\ref{s22},
and suppose that\/ $b$ is not an eigenvalue of\/ $a\De$ on $(M',g')$
in \S\ref{s21}, and on neither $(M',g')$ nor $(M'',g'')$ in \S\ref{s22}.
Then there exists $B>0$ such that for small\/ $t$, whenever $\xi\in
L^{2n/(n+2)}(M)$ there exists a unique $\phi\in L^2_1(M)$ with
$a\De\phi-b\phi=\xi$ on $(M,g_t)$,
and\/~${}^{g_t}\snm{\phi}{2}{1}\le B\,{}^{g_t}\lnm{\xi}{2n/(n+2)}$.
\label{l332}
\end{lem}

\begin{proof} The spectrum of $a\De$ on a compact manifold is discrete, 
so if $b$ is not an eigenvalue of $a\De$, then $a\De$ has no eigenvalues 
in a neighbourhood of $b$. Suppose $b$ is not an eigenvalue of 
$a\De$ on $(M',g')$ in \S\ref{s21}, and on neither $(M',g')$ nor
$(M'',g'')$ in \S\ref{s22}. Then there exists $\ga>0$ such that
$a\De$ has no eigenvalues in $(b-2\ga,b+2\ga)$ on these manifolds.
So by Theorem \ref{t331}, $a\De$ has no eigenvalues in $(b-\ga,b+\ga)$
on $(M,g_t)$ for small~$t$.

Thus easy analytical facts about the Laplacian imply that $a\De-b$ has a
right inverse $T:L^2(M)\ra L^2(M)$. As $M$ is compact, $L^2(M)\subset
L^{2n/(n+2)}(M)$. Let $\xi\in L^2(M)$. Then $\phi=T\xi\in L^2(M)$ exists and
satisfies the equation $a\De\phi-b\phi=\xi$ in the weak sense, and as
$a\De-b$ has no kernel, $\phi$ is unique. Since $\phi,\xi\in L^2(M)$,
multiplying this equation by $\phi$ and integrating gives a convergent
integral, so by subtraction, $\int_M\ms{\na\phi}\d V_{g_t}$ converges, and
$\phi\in L^2_1(M)$.

It remains to bound $\phi$ in $L^2_1(M)$. Let $\phi_1$ be the part of $\phi$
made up of eigenvectors of $a\De$ associated to eigenvalues less than $b$,
and $\phi_2$ the part associated to eigenvalues greater than $b$. Multiplying
the equation $a\De\phi-b\phi=\xi$ by $\phi_2-\phi_1$ and integrating gives
\e
\int_M\bigl(a\ms{\na\phi_2}-b\phi_2^2\bigr)\,\d V_{g_t}
-\int_M\bigl(a\ms{\na\phi_1}-b\phi_1^2\bigr)\,\d V_{g_t}
=\int_M(\phi_2-\phi_1)\xi\,\d V_{g_t}.
\label{s3eq10}
\e

But the restriction on the eigenvalues of $a\De$ means that
\begin{equation*}
\int_M\!a\ms{\na\phi_1}\d V_{g_t}\!\le\!(b\!-\!\ga)
\int_M\!\phi_1^2\d V_{g_t} \;\>\text{and}\;\>
\int_M\!a\ms{\na\phi_2}\d V_{g_t}\!\ge\!(b\!+\!\ga)
\int_M\!\phi_2^2\d V_{g_t},
\end{equation*}
and these together with \eq{s3eq10} and H\"older's inequality imply that
\begin{equation*}
\frac{\ga a}{a+b+\ga}\int_M\bigl(\ms{\na\phi_1}+
\ms{\na\phi_2}+\phi_1^2+\phi_2^2\bigr)\,\d V_{g_t}
\le{}^{g_t}\lnm{\phi_2-\phi_1}{p}{}^{g_t}\lnm{\xi}{2n/(n+2)}.
\end{equation*}

For $t\le\ze$, we apply Proposition \ref{p241} to $\phi_2-\phi_1$
to give ${}^{g_t}\lnm{\phi_2-\phi_1}{p} \le A\,{}^{g_t}\snm{\phi_2-
\phi_1}{2}{1}$. But $\phi_1,\phi_2$ are orthogonal in $L^2_1(M)$, so
${}^{g_t}\snm{\phi_2-\phi_1}{2}{1}={}^{g_t}\snm{\phi}{2}{1}$, and
similarly, the integral on the l.h.s.\ above is
${}^{g_t}\snm{\phi}{2}{1}^2$. Therefore
\begin{equation*}
\frac{\ga a}{a+b+\ga}{}^{g_t}\snm{\phi}{2}{1}^2\le
A\,{}^{g_t}\snm{\phi}{2}{1}\,{}^{g_t}\lnm{\xi}{2n/(n+2)}.
\end{equation*}

Dividing by $\ga a{}^{g_t}\snm{\phi}{2}{1}/(a+b+\ga)$ then gives
${}^{g_t}\snm{\phi}{2}{1}\le B\,{}^{g_t}\lnm{\xi}{2n/(n+2)}$ for small
$t$, where $B=(a+b+\ga)A/a\ga$. So the lemma holds for $\xi\in L^2(M)$.
This may be extended to $\xi\in L^{2n/(n+2)}(M)$ as in the proof of 
Lemma \ref{l321}, and the argument is complete.
\end{proof}

We now prove two existence theorems for metrics of scalar curvature 1:

\begin{thm} Let\/ $(M',g')$ and\/ $(M'',g'')$ be compact Riemannian
$n$-manifolds with scalar curvature $1$. Suppose that\/ $b$ is not an
eigenvalue of\/ $a\De$ on $(M',g')$, and that\/ $M',M''$ contain points
$m',m''$ with neighbourhoods in which $g',g''$ are conformally flat.

As in \S\ref{s21}, define the family $\{g_t:t\in (0,\de)\}$ of
metrics on $M=M'\# M''$. Then there exists $C>0$ such that\/ $g_t$
admits a smooth conformal rescaling to $\ti g_t=(1+\phi)^{p-2}g_t$
with scalar curvature $1$ for small\/ $t$,
and\/~${}^{g_t}\snm{\phi}{2}{1}\le Ct^2$.
\label{t333}
\end{thm}

\begin{thm} Let\/ $(M',g')$ and\/ $(M'',g'')$ be compact 
Riemannian $n$-manifolds with scalar curvature $1$. Suppose that\/
$b$ is not an eigenvalue of\/ $a\De$ on $(M',g')$ or $(M'',g'')$,
and that\/ $M',M''$ contain points $m',m''$ with neighbourhoods
in which $g',g''$ are conformally flat.

As in \S\ref{s22}, define the family $\{g_t:t\in (0,\de)\}$ of
metrics on $M=M'\# M''$. Then there exists $C>0$ such that\/ $g_t$
admits a smooth conformal rescaling to $\ti g_t=(1+\phi)^{p-2}g_t$
with scalar curvature $1$ for small\/ $t$,
and\/~${}^{g_t}\snm{\phi}{2}{1}\le Ct^2$.
\label{t334}
\end{thm}

\begin{proof}[Theorems \ref{t333} and \ref{t334}] These are the same as the 
proofs of Theorems \ref{t322} and \ref{t323}, except that Lemma \ref{l332}
should be applied in place of Lemma \ref{l321}, and where the proofs of
Theorems \ref{t322} and \ref{t323} mention scalar curvature $-1$, these
proofs should have scalar curvature~1.
\end{proof}

In \S\ref{s23} we saw that when $g',g''$ have scalar curvature 1,
\S\ref{s21} and \S\ref{s22} define three families of metrics $g_t$
on $M=M'\# M''$ with the same family of conformal classes $[g_t]$,
where $t$ in \S\ref{s22} corresponds to $t^2$ in \S\ref{s21}. Thus
Theorems \ref{t333} and \ref{t334} construct {\it three different}
metrics of scalar curvature 1 in the {\it same} conformal class
$[g_t]$ on $M$, for small $t$. The first resembles $(M',g')$ with
a small, asymptotically flat copy of $M''$ glued in at one point,
as in \S\ref{s21}. The second is like the first, but swapping
$M'$ and $M''$. The third resembles $(M',g')$ and $(M'',g'')$
joined by a small neck, as in~\S\ref{s22}.

These metrics are stationary points of the Hilbert action $Q$, but 
they need not be absolute minima, i.e.\ Yamabe metrics. The third 
is never a minimum. If $g'$ and $g''$ are Yamabe metrics, the author
expects that generically one of the first and second metrics will
be a Yamabe metric, and in a codimension 1 set of cases when
$\vol(M')\approx\vol(M'')$ both the first and second metrics will
be distinct Yamabe metrics.

\subsection{Extending to conformally curved metrics}
\label{s34}

In defining the metrics $\{g_t:t\in(0,\de)\}$ in \S\ref{s21} and
\S\ref{s22}, we assumed for simplicity that $g',g''$ are conformally
flat in neighbourhoods of $m',m''$. It turns out that if we drop this
assumption, then provided the metrics $g_t$ on $M$ are suitably defined
the results of \S\ref{s32} and \S\ref{s33} still hold without change.
The principal difference is that the expression \eq{s2eq2} for the
scalar curvature of $g_t$ becomes more complicated, with new error
terms that have to be estimated and controlled.

Following the method of \S\ref{s21}, we can choose an identification
of a ball about $m'$ in $M'$ with $B_{\de}(0)$, such that the induced
metric on $B_{\de}(0)$ is $g'=h+q'$, where $h$ is the standard metric
on $\R^n$ and $q'=O''(\md{v}^2)$ in the sense of \S\ref{s13}. To glue
in the neck metric $g_{\sst N}$ of \S\ref{s22}, for instance, we would
define $g_t=\be_1(h+q')+\be_2(1+t^{6(n-2)}\md{v}^{-(n-2)})^{p-2}h$,
where $(\be_1,\be_2)$ is the partition of unity defined in \S\ref{s21}.
Writing out the scalar curvature explicitly, we see that as 
$q'=O''(\md{v}^2)$, the terms involving $q'$ can be absorbed into
the existing error terms in \S\ref{s24}, so that \eq{s2eq3} holds
for the new metrics $g_t$. Therefore Proposition \ref{p231} holds
for the new metrics $g_t$ as well.

It can be seen by following the proof of Proposition \ref{p241} and the
proof of Theorem \ref{t331} in the appendix that no other nontrivial
modifications are required to prove these results for the more general
families of metrics $\{g_t:t\in(0,\de)\}$ discussed above. Thus the new
metrics satisfy all the necessary conditions, and the results of
\S\ref{s32} and \S\ref{s33} apply to them without change.

\section{Connected sums with zero scalar curvature}
\label{s4}

In this section the methods of \S\ref{s2} and \S\ref{s3} will be
adapted to study zero scalar curvature manifolds. We have three
cases to consider, when the scalar curvatures of $g'$ and $g''$ are
0 and 1, or 0 and 0, or $-1$ and 0. Each case introduces specific
difficulties, and each needs some additional methods to prove the 
existence of constant scalar curvature metrics.

The first two cases fit into a common analytic framework, and will be
handled together. The differences with the previous method are that 
$g_t$ must be defined more carefully than before to control the
errors sufficiently, and in the analysis the operator $a\De-\nu b$
now has one or two small eigenvalues. Thus when the sequence 
$\{\phi_i\}_{i=0}^\iy$ is defined inductively using the inverse of
this operator, the components in the directions of the corresponding
eigenvectors have to be attended to, to prevent the sequence diverging.
The third case is discussed in \S\ref{s45}. We shall outline what the
constant scalar curvature metrics look like and how to prove existence
results, but will not go into much detail.

\subsection{Combining zero and positive scalar curvature}
\label{s41}

Let $(M',g')$ and $(M'',g'')$ be compact Riemannian $n$-manifolds,
such that $g'$ has scalar curvature 0 and $g''$ scalar curvature 1.
Suppose as in \S\ref{s2} that $M',M''$ contain points $m',m''$
with neighbourhoods in which $g',g''$ are conformally flat.
As in \S\ref{s21} there exists a ball $B'$ about $m'$ in $M'$
and a diffeomorphism $\Phi':B_r(0)\subset\R^n\ra B'$ for $r<1$,
with $\Phi'(0)=m'$ and $(\Phi')^*(g')=(\psi')^{p-2}h$ for some
function $\psi'$ on $B_r(0)$ with $\psi'(0)=1$ and $\d\psi'(0)=0$.
As $g'$ has zero scalar curvature we have $\De\psi'=0$ by~\eq{s1eq2}.

Let $(\hat M,\hat g)$ be the stereographic projection of $M''$ from
$m''$. Then as in \S\ref{s21} there is an immersion $\Xi'':\R^n\sm
\ovB_R(0)\ra\hat M$ for some $R>0$, whose image is the complement
of a compact set in $\hat M$, such that $(\Xi'')^*(\hat g)=\xi^{p-2}h$,
where $\xi$ is a smooth function on $\R^n\sm\ovB_R(0)$ satisfying
$\xi(v)=1+O'(\md{v}^{2-n})$. Also, by Proposition \ref{p134} there
exists $\mu\in\R$ such that
\e
\xi(v)=1+\mu\md{v}^{2-n}+O'(\md{v}^{1-n}).
\label{s4eq1}
\e

As $(M'',g'')$ is not conformal to ${\cal S}^n$ with its round metric,
Theorem \ref{t135} shows that $\mu>0$, at least if $M''$ is spin or
$n\le 7$. For our purposes we only need to assume that $\mu\ne 0$,
so that $\mu\md{v}^{2-n}$ is the leading `error term' in
\eq{s4eq1}. Therefore if $n>7$ and $M''$ is not spin we suppose
that~$\mu\ne 0$.

Choose $k$ with $(n-2)(n+2)/2(n+1)<k<(n-2)(n+2)/2n$, which will
remain fixed throughout this section. Choose $\de\in(0,1)$ such
that $\de^{2k/(n+2)}\le r$ and $\de^{-2/n}\ge R$. We will write
down a family of metrics $\{g_t:t\in (0,\de)\}$ on $M=M'\# M''$,
in a similar way to \S\ref{s21}. For any $t\in (0,\de)$, define
$M$ and the conformal class of $g_t$ by
\begin{equation*}
M=\Bigl(M'\sm\Phi'\bigl[\ovB_{t^{(n-2)/n}}(0)\bigr]\Bigr)\amalg
\Bigl(\hat M\sm\bigl(\Xi''\bigl[\R^n\sm B_{t^{2k/(n+2)\!-\!1}}(0)
\bigr]\bigr)\Bigr)\Bigl/\bigti_t,
\end{equation*} 
where $\bigti_t$ is the equivalence relation defined by
\begin{equation*}
\Phi'\left[v\right]\,\,\bigti_t\,\,\Xi''[t^{-1}v]
\quad\text{whenever $v\in\R^n$ and $t^{(n-2)/n}<\md{v}<t^{2k/(n+2)}$.}
\end{equation*}

As in \S\ref{s21}, the conformal class $[g_t]$ of $g_t$ is the restriction 
of the conformal classes of $g'$ and $\hat g$ to the open sets of $M',\hat M$ 
making up $M$, and is well-defined because the conformal classes agree 
on the annulus of overlap $A_t$, where the two open sets are glued by 
$\bigti_t$. Define $g_t$ within this conformal class by $g_t=g'$ on 
the component of $M\sm A_t$ coming from $M'$, and $g_t=t^2\hat g$ on the 
component coming from $\hat M$. It remains to choose a conformal factor
on $A_t$. This is done as in \S\ref{s21}, except that the annulus
$\{v\in\R^n:t^{(n-2)/n}<\md{v}<t^{2k/(n+2)}\}$ in $\R^n$ replaces
$\{v\in\R^n:t^2<\md{v}<t\}$ in $\R^n$ in the definition of the partition
of unity. This completes the definition of $g_t$ for~$t\in(0,\de)$.

\begin{lem} Let the scalar curvature of\/ $g_t$ be $-\ep_t$. Then
$\ep_t$ is zero outside $A_t$. There exists $Y>0$ such that for all\/
$t\in(0,\de)$, $\ep_t$ satisfies $\md{\ep_t}\le Y$, and the volume of\/
$A_t$ with respect to $g_t$ satisfies $\vol(A_t)=O(t^{2nk/(n+2)})$.
Therefore ${}^{g_t}\lnm{\ep_t}{2n/(n+2)}=O(t^k)$
and\/~${}^{g_t}\lnm{\ep_t}{n/2}=O(t^{4k/(n+2)})$.
\label{l411}
\end{lem}

\begin{proof} Outside $A_t$, the metric $g_t$ is equal to $g'$ or 
homothetic to $\hat g$, and so has zero scalar curvature, verifying the 
first claim of the lemma. The proof that $\md{\ep_t}\le Y$ is the 
same as that for the corresponding statement in Proposition \ref{p231},
setting $\nu=0$. The estimate on the volume of $A_t$ also follows by
the method used in Proposition \ref{p231}, and the last two estimates
are immediate.
\end{proof}

We introduced $k$ above to make the estimate on
${}^{g_t}\lnm{\ep}{2n/(n+2)}$ easy to write down.

\begin{lem} For small\/ $t$,
\begin{equation*}
\int_M\ep_t\,\d V_{g_t}=(n-2)\om_{n-1}\mu t^{n-2}+O(t^{n-2+\al}),
\end{equation*}
where $\mu$ is the constant of\/ \eq{s4eq1}, $\om_{n-1}$ is
the volume of the unit sphere ${\cal S}^{n-1}$ in $\R^n$,
and\/~$\al=\min\bigl(2/n,2k(n+1)/(n+2)-(n-2)\bigr)>0$.
\label{l412}
\end{lem}

\begin{proof} Calculating with \eq{s1eq2} gives
\begin{align*}
\ep_t(v)=\psi_t^{-(n+2)/(n-2)}(v)\Bigl(&
2\bigl(\na\be_1(v)\bigr)\cdot
\bigl(\na\bigl(\psi'(v)-\xi(t^{-1}v)\bigr)\bigr)\\
&-\bigl(\De\be_1(v)\bigr)
\bigl(\psi'(v)-\xi(t^{-1}v)\bigr)\Bigr).
\end{align*}

Let $F$ be the quadratic form on $\R^n$ given by the second derivatives
of $\psi'$; then $\psi'=1+F+O'(\md{v}^3)$. As the scalar curvature of 
$g'$ is zero, the trace of $F$ is zero. Now 
$\d V_{g_t}=\psi_t^p\d V_h$. Multiplying through by this equation 
and making various estimates gives that
\e
\begin{split}
\ep_t(v)\,\d V_{g_t}=&\Bigl(2\bigl(\na\be_1(v)\bigr)\cdot
\bigl(\na F+O(\md{v}^2)-\mu t^{n-2}\na(\md{v}^{2-n})
-t^{n-1}O(\md{v}^{-n})\bigr)\\
&-\bigl(\De\be_1(v)\bigr)
\bigl(F+O(\md{v}^3)-\mu t^{n-2}\md{v}^{2-n}-t^{n-1}O(\md{v}^{1-n})
\bigr)\Bigr)\,\d V_h.
\end{split} 
\label{s4eq2}
\e

Integrate this over $A_t$. Now $\be_1(v)=\be(\md{v})$ where
$\be(\md{v})=\si(\log\md{v}/\log t)$, from \S\ref{s21}. So
$(\na\be_1(v))\cdot(\na F)=2\md{v}^{-1}F\frac{\d\be}{\d x}$ and
$(\na\be_1(v))\cdot(-\mu t^{n-2}\na(\md{v}^{2-n}))=(n-2)\mu t^{n-2}
\md{v}^{1-n}\frac{\d\be}{\d x}$ and $\De\be_1=-\frac{\d^2\be}{\d x^2}
+(1-n)\md{v}^{-1}\frac{\d\be}{\d x}$. Therefore
\e
\begin{split}
\int_{A_t}\ep_t(v)\,\d V_{g_t}=\int_{A_t}&\biggl(\frac{\d\be}{\d x}
\bigl((n+3)\md{v}^{-1}F+(n-3)\mu t^{n-2}\md{v}^{1-n}+O(\md{v}^2)\\
&+t^{n-1}O(\md{v}^{-n})\bigr)+\md{v}\frac{\d^2\be}{\d x^2}
\bigl(\md{v}^{-1}F-\mu t^{n-2}\md{v}^{1-n}\\
&\qquad\qquad\quad +O(\md{v}^3)+t^{n-1}O(\md{v}^{1-n})\bigr)\biggr)\,\d V_h.
\end{split}
\label{s4eq3}
\e

Using a Fubini theorem, we may write the integral on the right hand side as 
a double integral over ${\cal S}^{n-1}$ and $\md{v}$. The volume forms are
related by $\d V_h=\md{v}^{n-1}\d\Omega\d\md{v}$, where $\d\Omega$ is the 
standard volume form on ${\cal S}^{n-1}$ with radius 1. But as the trace of 
$F$ w.r.t.~$h$ is zero, $\int_{{\cal S}^{n-1}}F\d\Omega=0$ and the terms on 
the right hand side of \eq{s4eq3} involving $F$ vanish. So viewing 
\eq{s4eq3} as a double integral and integrating over ${\cal S}^{n-1}$ gives
\begin{equation*}
\int_{A_t}\!\!\ep_t(v)\d V_{g_t}\!=\!\om_{n-1}\int_{t^{(n\!-\!2)/n}
}^{t^{2k/(n+2)}}\!\Bigl((n\!-\!3)\mu t^{n\!-\!2}\frac{\d\be}{\d x}
\!-\!\mu t^{n\!-\!2}x\frac{\d^2\be}{\d x^2}\Bigr)\d x
\!+\text{error terms},
\end{equation*}
where $\om_{n-1}$ is the volume of ${\cal S}^{n-1}$.

The integral on the right is an exact integral, by parts. By definition, 
$\be$ changes from 0 to 1 and $\frac{\d\be}{\d x}$ from 0 to 0 over the 
interval, so it evaluates to
\begin{equation*}
\int_{A_t}\ep_t(v)\,\d V_{g_t}=(n-2)\om_{n-1}\mu t^{n-2}
+\text{\ error terms},
\end{equation*}
which is nearly the conclusion of the lemma; it remains only to 
show that the `error terms' are of order~$t^{n-2+\al}$. 

This is a simple calculation and will be left to the reader, the 
necessary ingredients being that as $\md{v}$ lies between $t^{(n-2)/n}$ and
$t^{2k/(n+2)}$, $O(\md{v})$ may be replaced by $O\bigl(t^{2k/(n+2)}\bigr)$ 
and $tO\bigl(\md{v}^{-1}\bigr)$ may be replaced by $O\bigl(t^{2/n}\bigr)$,
$\frac{\d\be}{\d x}=O\bigl(\md{v}^{-1}\bigr)$ and $\frac{\d^2\be}{\d x^2}=
O\bigl(\md{v}^{-2}\bigr)$. The error term that is usually the biggest is
$O\bigl(\md{v}^{n+1}\bigr)$, and in order to ensure that this error term is
smaller than the leading term calculated above, that is, to ensure $\al>0$,
$k$ must satisfy $k>(n+2)(n-2)/2(n+1)$, which was one of the conditions in 
the definition of $k$ above. 
\end{proof}

The lemma shows that the average scalar curvature of $(M,g_t)$ is 
close to $-(n-2)\om_{n-1}\mu t^{n-2}\vol(M')^{-1}$, so we will
choose this value for the scalar curvature of $\ti g_t$ in~\S\ref{s44}.

\subsection{Combining two metrics of zero scalar curvature}
\label{s42}

Let $M',M'',M,g',g'',m'$ and $m''$ be as usual, with $g',g''$ of 
zero scalar curvature and conformally flat in neighbourhoods of $m',m''$. 
Rescaling $g',g''$ by homotheties still gives metrics of zero scalar 
curvature. Thus gluing these rescaled metrics using the method of \S\ref{s22} 
gives a 2-parameter family of metrics in the same conformal class $[g_t]$ 
that all have small scalar curvature. Which of these metrics do we expect 
to be close to a metric of constant scalar curvature?

The necessary condition is that $\vol(M')=\vol(M'')$; we will explain
why at the end of \S\ref{s43}. Suppose, by applying a homothety to
$g'$ or $g''$ if necessary, that $\vol(M')=\vol(M'')$. A family of
metrics $\{g_t:t\in (0,\de)\}$ will be defined on $M$ following
\S\ref{s41}, such that when $t$ is small, $g_t$ resembles the union
of $M'$ and $M''$ with their metrics $g'$ and $g''$, joined by a small
`neck' of approximate radius $t$, which is modelled upon the manifold
$N$ of \S\ref{s22}, with metric~$t^2g_{\sst N}$.

Choose $k$ with $(n-2)(n+2)/2(n+1)<k<(n-2)(n+2)/2n$, and apply the
gluing method of \S\ref{s41} twice, once to glue one asymptotically
flat end of $(N,t^{2}g_{\sst N})$ into $M'$ at $m'$, and once to glue
the other asymptotically flat end into $M''$ at $m''$. The r\^ole of
$A_t$ in \S\ref{s41} is played by $A_t=A_t'\cup A_t''$, the disjoint
union of annuli $A_t'$ joining $N$ and $M'$ and $A_t''$ joining $N$
and $M''$. With this definition we state the next two lemmas,
which are analogues of Lemmas \ref{l411} and~\ref{l412}.

\begin{lem} Let the scalar curvature of\/ $g_t$ be $-\ep_t$. Then
$\ep_t$ is zero outside $A_t$. There exists $Y>0$ such that for all\/
$t\in (0,\de)$, $\ep_t$ satisfies $\md{\ep_t}\le Y$, and the volume
of\/ $A_t$ with respect to $g_t$ satisfies $\vol(A_t)=O(t^{2nk/(n+2)})$.
Therefore ${}^{g_t}\lnm{\ep_t}{2n/(n+2)}=O(t^k)$
and\/~${}^{g_t}\lnm{\ep_t}{n/2}=O(t^{4k/(n+2)})$.
\label{l421}
\end{lem}

\begin{proof} This is identical to Lemma \ref{l411}, and its proof is
the same, except that $g_t$ may also be homothetic to $g_{\sst N}$ 
in the first sentence. 
\end{proof}

\begin{lem} For all small\/ $t$,
\begin{align*}
&\int_{A_t'}\ep_t\,\d V_{g_t}=(n-2)\om_{n-1}t^{n-2}+O(t^{n-2+\al})\\
\text{and}\quad
&\int_{A_t''}\ep_t\,\d V_{g_t}=(n-2)\om_{n-1}t^{n-2}+O(t^{n-2+\al}),
\end{align*}
where $\al=\min\bigl(2/n,2k(n+1)/(n+2)-(n-2)\bigr)>0$ and\/ $\om_{n-1}$
is the volume of the unit sphere ${\cal S}^{n-1}$ in~$\R^n$.
\label{l422}
\end{lem}

\begin{proof} This is merely Lemma \ref{l412} applied twice, firstly
to the gluing of $N$ into $M'$ and secondly to the gluing of $N$ into
$M''$. We have also used the observation that for both asymptotically
flat ends of $N$, the constant $\mu$ of \S\ref{s41} takes the value 1.
To see this, compare the definition of $\mu$ in Proposition \ref{p134}
with the definition of $(N,g_{\sst N})$ in \S\ref{s22}.
\end{proof}

\subsection{Inequalities on the connected sum}
\label{s43}

Now we derive the analytic inequalities needed for the main result in
\S\ref{s44}. First of all, observe that Proposition \ref{p241} applies
to the metrics of \S\ref{s41} and \S\ref{s42}:

\begin{lem} Let\/ $\{g_t:t\in (0,\de)\}$ be one of the families of
metrics defined on $M=M'\# M''$ in \S\ref{s41} and\/ \S\ref{s42}.
Then there exist\/ $A>0$ and\/ $\ze\in(0,\de)$ such that\/
${}^{g_t}\lnm{\phi}{p}\le A\cdot{}^{g_t}\snm{\phi}{2}{1}$
whenever $\phi\in L^2_1(M)$ and\/~$t\in(0,\ze]$.
\label{l431}
\end{lem}

\begin{proof} The proof follows that of Proposition \ref{p241}, applied to the
metrics of \S\ref{s41} and \S\ref{s42} rather than \S\ref{s21} and \S\ref{s22}, except for
some simple changes to take into account the different powers of $t$ 
used to define the new metrics. 
\end{proof}

As in \S\ref{s33}, to calculate with the inverse of $a\De-\nu b$ we need to 
know about the spectrum of $a\De$ on $(M,g_t)$. The next three results 
give the necessary information; the proofs are again deferred to the
appendix.

\begin{thm} Let\/ $\{g_t:t\in (0,\de)\}$ be the family of metrics
defined on $M=M'\# M''$ in \S\ref{s41}. Choose $\ga>0$ such that\/
$a\De$ on $(M',g')$ has no eigenvalues in $(0,2\ga)$. Then $a\De$
on $(M,g_t)$ has no eigenvalues in $(0,\ga)$ for small\/~$t$.
\label{t432}
\end{thm}

For the metrics $g_t$ of \S\ref{s42} the situation is more complicated.
For small $t$ we expect the eigenvectors of $a\De$ on $(M,g_t)$ with
small eigenvalues to be close to eigenvectors of $a\De$ on $M'$ or
$M''$ with small eigenvalues, that is, to constant functions on $M'$
and $M''$. So we expect two eigenvectors on $(M,g_t)$ associated to
small eigenvalues, one the constant function, and the other close to
a constant on the $M'$ part of $M$, and to a different constant on
the $M''$ part of $M$. The next result describes this second
eigenvector. The proof is deferred to the appendix.

\begin{prop} Let\/ $\{g_t:t\in (0,\de)\}$ be the family of 
metrics defined on $M=M'\# M''$ in \S\ref{s42}. Then for
small\/ $t$ there exists $\la_t>0$ and\/ $\be_t\in C^\iy(M)$
such that\/ $a\De\be_t=\la_t\be_t$ on $(M,g_t)$. Here
$\la_t=O(t^{n-2})$, and
\begin{equation*}
\be_t=\begin{cases} 
\phantom{-}1+O(t^{n-2}) &\text{on $M'\sm B'$}\\
\phantom{-}1+O(t^{n-2}\md{v}^{2-n}) 
&\text{on $\{v:t\le\md{v}<\de\}\subset B'$}\\
-1+O(t^{n-2}) & \text{on $M''\sm B''$}\\
-1+O(t^{n-2}\md{v}^{2-n}) 
&\text{on $\{v:t\le\md{v}<\de\}\subset B''$,}
\end{cases}
\end{equation*}
identifying subsets of\/ $M',M''$ with subsets of\/~$M$.
\label{p433}
\end{prop}

The proposition is proved by a series method, starting with a function 
that is 1 on the part of $M$ coming from $M'$ and $-1$ on the part coming 
from $M''$, and then adding small corrections to get to an eigenvector
of $a\De$. Note that $\be_t$ takes the approximate values $\pm 1$ on the
two halves because $\vol(M')=\vol(M'')$ by assumption; if the volumes 
were different, then the approximate values would have to be adjusted 
so that~$\int_M\be_t\d V_{g_t}=0$.

We may now state the analogue of Theorem \ref{t432} for the metrics of
\S\ref{s42}, which will be proved in the appendix.

\begin{thm} Let\/ $\{g_t:t\in (0,\de)\}$ be the family of metrics
defined on $M=M'\# M''$ in \S\ref{s42}. Choose $\ga>0$ such that all
$a\De$ has no eigenvalues in $(0,2\ga)$ on $(M',g')$ or $(M'',g'')$.
Then for small\/ $t$, the only eigenvalue of $a\De$ in $(0,\ga)$ on
$(M,g_t)$ is $\la_t$ from Proposition \ref{p433}, with
eigenspace~$\an{\be_t}$.
\label{t434}
\end{thm}

Theorems \ref{t432} and \ref{t434} will fit into the existence proofs
of \S\ref{s44} in the same way as Theorem \ref{t331} does into that of
\S\ref{s33}. The small eigenvalue $\la_t$ in Theorem \ref{t434} means
that $\be_t$-components of functions will have to be carefully controlled,
to ensure that inverting $a\De$ upon them does not give a large result.
We now bound the $\be_t$-component of~$\ep_t$.

\begin{lem} Let\/ $\{g_t:t\in (0,\de)\}$ be the family of metrics
defined on $M=M'\# M''$ in \S\ref{s42}. Then for small\/~$t$,
\begin{equation*}
\int_M\be_t\ep_t\,\d V_{g_t}=O(t^{n-2+\al}),
\end{equation*}
where $\be_t$ is the function of Proposition \ref{p433} and\/
$\al$ is as in Lemma~\ref{l422}.
\label{l435}
\end{lem}

\begin{proof} Proposition \ref{p433} shows that $\be_t=1+O(t^{2(n-2)/n})$
on $A'_t$ and $\be_t=-1+O(t^{2(n-2)/n})$ on $A''_t$, as these are annuli
in which $t^{(n-2)/n}<\md{v}<t^{2k/(n+2)}$. Applying these and Lemmas
\ref{l421} and \ref{l422} to the integral of $\be_t\ep_t$ over $M$ gives
\begin{equation*}
\int_M\be_t\ep_t\,\d V_{g_t}=O(t^{n-2+\al})
+Y\vol(A_t)\cdot O(t^{2(n-2)/n}),
\end{equation*}
and as $\vol(A_t)=O(t^{2nk/(n+2)})$, the second term is
$O(t^{2nk/(n+2)+2(n-2)/n})$. But by the definitions of $k$ and $\al$, 
it is easily shown that $n-2+\al<2nk/(n+2)+2(n-2)/n$, and so the first 
error term is larger and subsumes the second, as required. 
\end{proof}

We note that this lemma is the reason for requiring that
$\vol(M')=\vol(M'')$. For if the two are not equal, then
Lemma \ref{l422} still shows that $\int_{A'_t}\ep_t\d V_{g_t}$ and
$\int_{A''_t}\ep_t\d V_{g_t}$ are equal to highest order, but
$\be_t$ takes values approximately proportional to $\vol(M')^{-1}$ on
$A'_t$, and to $\vol(M'')^{-1}$ on $A''_t$. Thus in this case
$\int_M\be_t\ep_t\d V_{g_t}$ is $O(t^{n-2})$ rather than
$O(t^{n-2+\al})$. But we will need $\int_M\be_t\ep_t\d V_{g_t}
=o(t^{n-2})$ for the proof in~\S\ref{s44}.

\subsection{Existence of constant scalar curvature metrics}
\label{s44}

Now we give the existence results for constant scalar curvature metrics
on the connected sums of \S\ref{s41} and~\S\ref{s42}.

\begin{thm} Let\/ $\{g_t:t\in (0,\de)\}$ be one of the families of
metrics on $M=M'\# M''$ defined in \S\ref{s41} or \S\ref{s42}. Then
there exists $C>0$ such that\/ $g_t$ admits a smooth conformal rescaling
for small\/ $t$ to $\ti g_t=(1+\phi)^{p-2}g_t$, with scalar curvature
$-(n\!-\!2)\om_{n-1}\mu t^{n-2}\vol(M')^{-1}$ in \S\ref{s41} and\/
$-(n\!-\!2)\om_{n-1}t^{n-2}\vol(M')^{-1}$ in \S\ref{s42}, and\/
${}^{g_t}\snm{\phi}{2}{1}\le Ct^{\al}$. Here $\om_{n-1}$ is the
volume of the unit sphere ${\cal S}^{n-1}$ in $\R^n$, and\/ $\mu,\al$
are as in \S\ref{s41}, and we suppose $\mu\ne 0$ in the cases not
covered by Theorem~\ref{t135}.
\label{t441}
\end{thm}

\begin{proof} Let $D_0$ be $(n-2)\om_{n-1}\mu\vol(M')^{-1}$ in \S\ref{s41}
and $(n-2)\om_{n-1}\vol(M')^{-1}$ in \S\ref{s42}. Define a function $\eta$
on $(M,g_t)$ by $\eta=\ep_t-D_0t^{n-2}$. Then $g_t$ has scalar curvature
$-D_0t^{n-2}-\eta$. As in \S\ref{s3}, the condition for $\ti g_t=(1+\rho+
\tau)^{p-2}g_t$ to have scalar curvature $-D_0t^{n-2}$ is
\e
a\De(\rho+\tau)+bD_0t^{n-2}(\rho+\tau)=\eta+\eta\cdot(\rho+\tau)
-D_0t^{n-2}f(\rho+\tau).
\label{s4eq4}
\e

Define a vector space $P$ of functions on $(M,g_t)$ by
$P=\an{1}$ in \S\ref{s41}, and $P=\an{1,\be_t}$ in \S\ref{s42},
where $\be_t$ is as in Proposition \ref{p433}. We shall construct
$\rho,\tau$ satisfying \eq{s4eq4}, with $\rho\in P$ and $\tau\in
P^\perp$ with respect to the $L^2_1$ inner product.

Define inductively sequences $\{\rho_i\}_{i=0}^\iy$ of elements of $P$
and $\{\tau_i\}_{i=0}^\iy$ of elements of $P^\perp\subset L^2_1(M)$ by
$\rho_0=\tau_0=0$, and having defined the sequences up to $i-1$, let $\rho_i$
and $\tau_i$ be the unique elements of $P$ and $P^\perp$ satisfying
\e
\begin{split}
a\De(\rho_i+\tau_i)+bD_0t^{n-2}(\rho_i+\tau_i)=\,&\eta+\eta\cdot
(\rho_{i-1}+\tau_{i-1})\\
&-D_0t^{n-2}f(\rho_{i-1}+\tau_{i-1}).
\end{split}
\label{s4eq5}
\e
If we can show that these sequences converge to $\rho\in P$ and
$\tau\in P^\perp$ that are small when $t$ is small, then the arguments 
of \S\ref{s3} complete the theorem.

The difficulty lies in inverting the operator $a\De+bD_0t^{n-2}$: by
Theorems \ref{t432} and \ref{t434}, the operator is invertible on
$P^\perp$ with inverse bounded by $\ga^{-1}$, as all the eigenvectors
of $a\De$ in $P^\perp$ have eigenvalues at least $\ga$. But on $P$,
the inverse is of order $t^{2-n}$, which is large; so $\rho_i$ may be
large even if the right hand side of \eq{s4eq5} is small.

The solution is to ensure that the $P$ components of $\eta$ are smaller even
than $t^{n-2}$, so that after applying the inverse of $a\De+bD_0t^{n-2}$ to
them, they are still small. Let $\pi$ denote orthogonal projection onto $P$;
both the $L^2$ and the $L^2_1$ inner product give the same answer, and in fact
the projection makes sense even in $L^1(M)$. Then from \eq{s4eq5} we make the
estimates 
\ea
\begin{split}
\snm{\rho_i}{2}{1}&\le D_1t^{2-n}\Bigl(\lnm{\pi(\eta)}{1}
+\lnm{\pi(\eta\rho_{i-1})}{1}
+\lnm{\pi(\eta\tau_{i-1})}{1}\Bigr)\\
&\qquad\qquad\qquad +D_2\blnm{\pi\bigl(f(\rho_{i-1}+\tau_{i-1})\bigr)}{1},
\end{split}
\label{s4eq6}\\
\begin{split}
\snm{\tau_i}{2}{1}&\le D_3\Bigl(\lnm{\eta}{2n/(n+2)}
+\lnm{\eta\rho_{i-1}}{2n/(n+2)}
+\lnm{\eta\tau_{i-1}}{2n/(n+2)}\\
&\qquad\qquad\qquad 
+D_0t^{n-2}\lnm{f(\rho_{i-1}+\tau_{i-1})}{2n/(n+2)}\Bigr),
\end{split}
\label{s4eq7}
\ea
for some constants $D_1,D_2,D_3$ independent of $t$. The norms on the 
right hand side of \eq{s4eq6} would normally be $L^{2n/(n+2)}$ norms, 
but as $P$ is a finite-dimensional space all norms are equivalent, and 
we may use the $L^1$ norm.

Our strategy is to show that if $\snm{\rho_{i-1}}{2}{1}\le D_4t^{\al}$ and
$\snm{\tau_{i-1}}{2}{1}\le D_5t^k$ for large enough constants $D_4,D_5$, then
$\snm{\rho_i}{2}{1}\le D_4t^{\al}$ and $\snm{\tau_i}{2}{1}\le D_5t^k$ also
hold for small $t$, so by induction the sequences are bounded; convergence
for small $t$ easily follows by a similar argument to that used in
Lemma~\ref{l312}.

From Lemmas \ref{l412}, \ref{l422} and \ref{l435} we deduce that
$\lnm{\pi(\eta)}{1}=O(t^{n-2+\al})$, so the first term on the right of
\eq{s4eq6} contributes $O(t^{\al})$ to $\snm{\rho_i}{2}{1}$, consistent with
$\snm{\rho_i}{2}{1}\le D_4t^{\al}$ if $D_4$ is chosen large enough. The
third term $\lnm{\pi(\eta\tau_{i-1})}{1}$ is bounded by
$A\lnm{\eta}{2n/(n+2)}\snm{\tau_{i-1}}{2}{1}$, and 
$\lnm{\eta}{2n/(n+2)}=O(t^k)$ by Lemmas \ref{l411} and \ref{l421}; the third term 
therefore contributes $O(t^{2k+2-n})$ to $\snm{\rho_i}{2}{1}$, and by the 
definition of $\al$, this error term is strictly smaller than $O(t^{\al})$. 
The fourth error term is also easily shown to be smaller than~$O(t^{\al})$.

Thus the only problem term in \eq{s4eq6} is the second term, and the only 
reason it is a problem is that the $P^\perp$ component of $\eta$, 
multiplied by $\rho_{i-1}$, may have an appreciable component in $P$. 
We get round this as follows. Suppose $\xi\in P^\perp$ and $\rho\in P$, 
and consider the $P$-component of $\xi\rho$. In \S\ref{s41} this 
component is zero, and there is no problem; in \S\ref{s42} there 
may be a component in the direction of $\be_t$, and it is measured by
$\int_M\xi\be_t^2\d V_{g_t}$. But by the description of $\be_t$ in
Proposition \ref{p433}, $\be_t^2$ is close to 1, and $\xi$ is
orthogonal to the constants, and so in general the $P$ component of
$\xi\rho$ will be small compared to the sizes of $\xi$ and $\rho$.
Taking this into account, it is easy to get a good bound
on~$\lnm{\pi(\eta\rho_{i-1})}{1}$.

The rest of the proof will be left to the reader. What remains to be done 
is to prove inductively that bounds $\snm{\rho_i}{2}{1}\le D_4t^{\al}$,
$\snm{\tau_i}{2}{1}\le D_5t^k$ hold for small enough $t$, and then to prove 
the convergence of the sequences, and these may both be done using the 
methods of Lemma \ref{l312}, working from \eq{s4eq6} and \eq{s4eq7}. Setting 
$\phi=\rho+\tau$, where $\rho,\tau$ are the limits of the sequences, the 
reader may then rejoin the proof of Theorem \ref{t311} after Lemma \ref{l313}. 
\end{proof}

As the metrics constructed have negative scalar curvature, they are 
unique in their conformal classes, and are Yamabe metrics. The theorem 
thus tells us that the Yamabe metric on the connected sum, with small 
neck, of two zero scalar curvature manifolds, balances the volumes of 
the two component manifolds so that they are equal, a fact which seems 
rather appealing.

\subsection{Combining zero and negative scalar curvature}
\label{s45}

There is just one case left, that of gluing a metric of zero scalar curvature 
into a metric of scalar curvature $-1$. This case can be handled using the
results of \S\ref{s3}, using the following simple extension of the method: we
define a family of metrics $\{g_t:t\in(0,\de)\}$ on the connected sum $M$,
with $-1-\ep_t$ the scalar curvature of $g_t$, and then prove that
$\md{\ep_t}\le Y$, $\lnm{\ep_t}{n/2}\le Zt^{\iota}$ as in Proposition
\ref{p231}, and that $\lnm{\phi}{p}\le At^{-\kappa}\snm{\phi}{2}{1}$ for
$\phi\in L^2_1(M)$ as in Proposition \ref{p241}, where norms are taken
with respect to~$g_t$. 

Here $\iota,\kappa>0$. The idea is that if $\iota$ is large compared to
$\kappa$, then we may follow the proofs of \S\ref{s3} adding in powers
of $t$, and at the crucial stages when we need some expression to be
sufficiently small, the power of $t$ will turn out to be positive,
and so we need only take $t$ small enough. To do this in practice,
we modify the proof slightly to cut down the number of applications
of Proposition \ref{p241}, and thus obtain a more favourable necessary
ratio of $\iota$ to~$\kappa$.

Now $\kappa$ is essentially determined by what the Yamabe metric on
the connected sum actually looks like -- if $t$ is the radius of the
`neck', any family of good approximations to the Yamabe metrics will
have $\kappa=(n-2)/n$ (this value will be justified below). So the
problem is to make $\iota$ large enough, in other words, to start
with a family of metrics $g_t$ that are a good approximation to
scalar curvature $-1$. We shall not go through the construction and
proof again for this case, but will describe how to define metrics
$g_t$ that have close enough to constant scalar curvature for the
method outlined above to work.

Consider the connected sum $M$ of $(M',g')$ with scalar curvature
$-1$ and $(M'',g'')$ with scalar curvature 0. To get a good enough
approximation to scalar curvature $-1$, we have to rescale $g''$
so that its scalar curvature approximates $-1$ rather than 0. Let
$\xi$ be the Green's function of $a\De$ at $m''$ on $M''$ satisfying
$a\De\xi=\de_{m''}-\vol(M'')^{-1}$ in the sense of distributions.
Since $\xi$ is only defined up to the addition of a constant, choose
$\xi$ to have minimum value 0. Then $\xi$ is a $C^\iy$ function on
$\hat M=M''\sm\{m''\}$ with a pole at $m''$, of the form
$(n-2)\om^{-1}_{n-1}\md{v}^{2-n}+O'(\md{v}^{1-n})$, in the usual
coordinates.

Let $\hat g_t=\bigl(t^{(n-2)^2/2n}+t^{(n-2)(n+2)/2n}\vol(M'')\xi
\bigr)^{p-2}g''$ on $\hat M$. Calculating its scalar curvature
$\hat S_t$ using \eq{s1eq2} gives
\begin{equation*}
\hat S_t=-\bigl(1+t^{2(n-2)/n}\vol(M'')\xi\bigr)^{-(n+2)/(n-2)},
\end{equation*}
so that $-1\le\hat S_t<0$, and $\hat S_t$ is close to $-1$ away from
$m''$ for small $t$. Outside a small neighbourhood of $m''$, $\hat g_t$
is close to $t^{2(n-2)/n}g''$, so that the diameter of $M''$ is multiplied
by $t^{(n-2)/n}$. But in a small neighbourhood of $m''$, $\hat g_t$
resembles a `neck' metric of radius proportional to $t$, as in
\S\ref{s22}. So $\hat g_t$ looks like $M''$ rescaled by $t^{(n-2)/n}$,
and with a `neck' of radius proportional to $t$, opening out to an
asymptotically flat end.

We construct $g_t$ by gluing $\hat g_t$ into $g'$ using the natural `neck'.
Thus a rough description of $g_t$ is that it is $g'$ on the $M'$ part 
and $t^{2(n-2)/n}g''$ on the $M''$ part, and the two parts are joined 
by a `neck' with radius proportional to $t$. With a family 
$\{g_t:t\in(0,\de)\}$ of metrics defined in this way, the modified 
method outlined above may be applied to show that there exist small 
conformal deformations of $g_t$ to scalar curvature $-1$, for $t$ 
sufficiently small.

One final point: we can now see the reason for the failure of Proposition
\ref{p241}, which necessitated this whole detour. Consider a smooth
function $\phi$ on $(M,g_t)$ that is 0 on the $M'$ part, 1 on the $M''$
part, and changes only on the `neck'. A simple calculation shows that
for such a function, $\lnm{\phi}{p}\sim t^{(2-n)/n}\snm{\phi}{2}{1}$,
and so the value for $\kappa$ given above is the least possible.

\subsection{Doing without conformal flatness}
\label{s46}

In \S\ref{s34}, we explained that the results of \S\ref{s3} still hold
if the assumption that $g',g''$ are conformally flat in neighbourhoods of 
$m',m''$ is dropped. However, the results of \S\ref{s4} do require 
modifications to generalize in this way. The problem is in extending 
Lemmas \ref{l412} and \ref{l422} to the curved case: we need a quite precise 
evaluation of the total scalar curvature of $g_t$, and have to be 
careful that the error terms do not swamp the term that we can evaluate.

To deal with the case of \S\ref{s41} first, it can be shown that the proof
of Lemma \ref{l412} still holds when $g'$ has conformal curvature near $m'$, 
because defining $g_t$ as in \S\ref{s34}, the extra error terms introduced
in the scalar curvature can be absorbed into the error terms already in 
\eq{s4eq2}, and so the proof of Lemma \ref{l412} holds from that point.
But if we allow $g''$ to be conformally curved near $m''$, then
Proposition \ref{p134} doesn't hold, and the scalar curvature of
$g_t$ is dominated (except for $n=3,4,5$) by error terms that seem
to have no nice expression.

Therefore the situation is this. Theorem \ref{t441} applies without change
when the metrics $g_t$ of \S\ref{s41} are defined using $g''$ conformally
flat near $m''$, but $g'$ not necessarily conformally flat near $m'$. To
include the case $g''$ not conformally flat, the result will hold if we
weaken it so as not to prescribe the constant value that the scalar
curvature takes, but merely give an estimate of its magnitude. Also,
I believe that the result applies as stated when $n=3,4$ or 5, because
then the mass term is large enough to dominate the errors.

The case of \S\ref{s42} is easier: Theorem \ref{t441} applies without
change to the metrics $g_t$ of \S\ref{s42} defined using $g', g''$ not
supposed conformally flat about $m'$ and $m''$. This is because the
metrics $g_t$ are made by gluing in the neck metric $t^2g_{\sst N}$,
which is conformally flat, so it reduces to the case of \S\ref{s41}
when $g''$ is locally conformally flat, which we have already seen works.

This leaves the case of \S\ref{s45}. The generalization of Proposition
\ref{p241} suggested there will extend without change to the conformally
curved case, so the problem is to define the metrics $g_t$ in such a way
that the extra error terms introduced in the expression for the scalar
curvature still give good enough approximations to constant scalar
curvature for the existence result to apply. I think this can be done
quite readily just by working on how to produce good approximations
$g_t$, say by adding a well-chosen conformal factor to the existing
definition.

\appendix
\section{Appendix. The spectrum of $a\De$ on $M'\# M''$}
\label{A}

In this appendix we prove Theorems \ref{t331}, \ref{t432} and \ref{t434},
and Proposition \ref{p433}. They are results on the eigenvalues and
eigenvectors of the operator $a\De$ on $M$ with the metrics $g_t$ defined
in \S\ref{s3} and \S\ref{s4}. They appear here and not in the main text
because the proofs are long calculations.

Theorem \ref{t331} takes up \S\ref{A1} and \S\ref{A2}. Its proof divides 
naturally into considering eigenvalues of $a\De$ smaller than $b$
and eigenvalues larger than $b$. The eigenvectors with eigenvalues smaller
than $b$ form a finite dimensional space $E$. In \S\ref{A1} we define a
vector space $E_t$ that is a good approximation to $E$ when $t$ is small,
and using this we show that all eigenvalues of vectors in $E$ are at most
$b-\ga$. In \S\ref{A2} we consider eigenvectors with eigenvalues larger
than $b$, which must therefore be orthogonal to $E$, and by considering
their inner product with $E_t$ we can show that their eigenvalues must be
at least $b+\ga$.

In \S\ref{A3}, we prove similar results for use in the zero
scalar curvature material of \S\ref{s4}. Most of the work needed to prove them
has already been done in \S\ref{A1} and \S\ref{A2}, and the main problem is
the construction of an eigenvector $\be_t$ of $a\De$ with a small
eigenvalue $\la_t$. This is done by a sequence method, the basic idea being
to start with an approximation to $\be_t$ and repeatedly invert $a\De$ upon
it; as $\la_t$ is the smallest positive eigenvalue, the $\be_t$-component
of the resulting sequence grows much faster than any other and so comes to
dominate.

\subsection{Eigenvalues of $a\De$ smaller than $b$}
\label{A1}  

We now prove Theorem \ref{t331}, which is reproduced here.

\begin{result}{Theorem \ref{t331}} Let\/ $\{g_t:t\in(0,\de)\}$ be one
of the families of metrics defined on $M=M'\# M''$ in \S\ref{s21} or
\S\ref{s22}, and suppose that for some $\ga>0$, $a\De$ has no eigenvalues
in $(b-2\ga,b+2\ga)$ on $(M',g')$ in \S\ref{s21}, and on both\/ $(M',g')$
and\/ $(M'',g'')$ in \S\ref{s22}. Then $a\De$ has no eigenvalues in
$(b-\ga,b+\ga)$ on $(M,g_t)$ for small\/~$t$.
\end{result}

\begin{proof} It is well known that the spectrum of the Laplacian on
a compact Riemannian manifold is discrete and nonnegative, and that 
the eigenspaces are finite-dimensional. Therefore on $M'$ and $M''$ 
there are only finitely many eigenvalues of $a\De$ smaller than $b$, 
and to each is associated a finite-dimensional space of eigenfunctions. 

Let $E'$ be the finite-dimensional vector space of smooth functions 
on $M'$ generated by eigenfunctions of $a\De$ on $M'$ associated to 
eigenvalues less than $b$; we think of $E'$ as a subspace of $L^2_1(M')$. 
For \S\ref{s22}, define $E''$ on $M''$ in the same way. As $a\De$ has
no eigenvalues in the interval $(b-2\ga,b+2\ga)$, we see that:
\ea
&\text{if $\phi\in E'$, then\ }
\int_{M'}a\ms{\na\phi}\,\d V_{g'} \le 
\int_{M'}(b-2\ga)\phi^2\,\d V_{g'},
\label{Aeq1}\\
&\text{if $\phi\in (E')^\perp\subset L^2_1(M')$, then\ }
\int_{M'}a\ms{\na\phi}\,\d V_{g'} \ge 
\int_{M'}(b+2\ga)\phi^2\,\d V_{g'},
\label{Aeq2}
\ea
and also two analogous inequalities for $M''$ in \S\ref{s22}. 
The perpendicular subspace $(E')^\perp$ of \eq{Aeq2} may be taken 
with respect to the inner product of $L^2_1(M')$ or with respect to 
that of $L^2(M')$ -- both give the same space, as $E'$ is a sum of 
eigenspaces of~$a\De$.

Now if we have two statements like \eq{Aeq1} and \eq{Aeq2} but applying 
to $M$ rather than $M'$, then we can prove the result. This is the
content of the next lemma.

\begin{lem} Suppose there is a subspace $E_t$ of\/ $L^2_1(M)$
for small $t$ such that
\ea
&\text{if\/ $\phi\in E_t$, then\ }
\int_Ma\ms{\na\phi}\,\d V_{g_t}\le 
\int_M(b-\ga)\phi^2\,\d V_{g_t},\quad\text{and}
\label{Aeq3}\\
&\text{if\/ $\phi\in (E_t)^\perp\subset L^2_1(M)$, then\ }
\int_Ma\ms{\na\phi}\,\d V_{g_t}\ge 
\int_M(b+\ga)\phi^2\,\d V_{g_t},
\label{Aeq4}
\ea
defining $(E_t)^\perp$ with the $L^2_1$ inner product. Then
Theorem \ref{t331} holds.
\label{lA11}
\end{lem}

\begin{proof} Suppose $E_t$ exists. We must show that $l\in
(b-\ga,b+\ga)$ cannot be an eigenvalue of $a\De$ on $(M,g_t)$.
Suppose $\phi$ is an eigenfunction of $a\De$ with this eigenvalue
$l$. Let $\phi_1$ and $\phi_2$ be the components of $\phi$ in
$E_t$ and $(E_t)^\perp$ respectively. Then, as $a\De\phi-l\phi=0$,
\begin{align*}
0&=\int_M(\phi_2-\phi_1)\bigl(a\De(\phi_1+\phi_2)
-l(\phi_1+\phi_2)\bigr)\,\d V_{g_t}\\
&=\int_M\bigl(a\ms{\na\phi_2}-l\phi_2^2\bigr)\,\d V_{g_t}
-\int_M\bigl(a\ms{\na\phi_1}-l\phi_1^2\bigr)\,\d V_{g_t}\\
&\ge \int_M\bigl((\ga+l-b)\phi_1^2+(\ga+b-l)\phi_2^2\bigr)
\,\d V_{g_t},
\end{align*}
using \eq{Aeq3} and \eq{Aeq4} in the last line. But as
$\ga+l-b,\ga+b-l>0$, this shows that $\phi_1=\phi_2=\phi=0$,
which is a contradiction. 
\end{proof}

To complete the proof of the theorem, we therefore need to produce some 
spaces $E_t$ of functions on $M$ satisfying \eq{Aeq3} and \eq{Aeq4}. In
\S\ref{s21}, $E_t$ should be modelled on $E'$, and in \S\ref{s22},
on~$E'\oplus E''$.

As a half-way stage between $E', E''$ and $E_t$, spaces $\ti E',\ti E''$ 
of functions on $M',M''$ will be made that are close to $E',E''$, but 
which vanish on small balls around $m',m''$. Let $\si'$ be a $C^\iy$
function on $M'$ that is 1 on the complement of a small ball about
$m'$, 0 on a smaller ball about $m'$, and otherwise taking values
in $[0,1]$. Now define~$\ti E'=\si'E'=\{\si'v:v\in E'\}$.

By choosing the ball outside which $\si'$ is 1 to be small, we can
ensure that $\ti E'$ is close to $E'$ in $L^2_1(M')$ in the following
sense: the two have the same dimension, and any $\ti v\in\ti E'$ may
be written as $\ti v=v_1+v_2$, where $v_1,v_2\in E',(E')^\perp$
respectively, and satisfy
\e
\snm{v_2}{2}{1}^2\le\frac{\ga}{2(a+b+2\ga)}\snm{v_1}{2}{1}^2.
\label{Aeq5}
\e

Suppose that $\si'$ has been chosen so that these hold. Then two
statements similar to \eq{Aeq1} and \eq{Aeq2} hold for $\ti E'$, as
we shall see in the next lemma.

\begin{lem} The subspace $\ti E'$ satisfies the following 
two conditions:
\ea
&\text{if\/ $\phi\in\ti E'$, then\ }
\int_{M'}a\ms{\na\phi}\,\d V_{g'} \le
\int_{M'}(b-{\ts\frac{3}{2}}\ga)\phi^2\,\d V_{g'},
\label{Aeq6}\\
&\text{if\/ $\phi\in(\ti E')^\perp\subset L^2_1(M')$,
then\ } \int_{M'}a\ms{\na\phi}\,\d V_{g'} \ge
\int_{M'}(b+{\ts\frac{3}{2}}\ga)\phi^2\,\d V_{g'},
\label{Aeq7}
\ea
where the inner product used to construct\/ $(\ti E')^\perp$ is that 
of\/~$L^2_1(M')$.
\label{lA12}
\end{lem}

\begin{proof} First we prove \eq{Aeq6}. Let $\phi\in\ti E'$; then 
$\phi=v_1+v_2$, with $v_1\in E'$ and $v_2\in(E')^\perp$. Because 
$v_1$ and $v_2$ are orthogonal in both $L^2$ and $L^2_1$,
\begin{align*}
a\int_{M'}&\ms{\na\phi}\,\d V_{g'}\\
&=a\int_{M'}\ms{\na v_1}\,\d V_{g'}+
a\int_{M'}\ms{\na v_2}\,\d V_{g'}\\ 
&\le a\int_{M'}\ms{\na v_1}\,\d V_{g'}+a\snm{v_2}{2}{1}^2\\ 
&\le a\int_{M'}\ms{\na v_1}\,\d V_{g'}+
\frac{a\ga}{2(a+b+2\ga)}\snm{v_1}{2}{1}^2\\
&=a\left(1+\frac{\ga}{2(a+b+2\ga)}\right)
\int_{M'}\ms{\na v_1}\,\d V_{g'}+\frac{a\ga}{2(a+b+2\ga)}
\int_{M'}v_1^2\,\d V_{g'}\\
&\le \left((b-2\ga)\left(1+\frac{\ga}{2(a+b+2\ga)}\right)
+\frac{a\ga}{2(a+b+2\ga)}\right)
\int_{M'}v_1^2\,\d V_{g'}\\ 
&\le (b-{\ts\frac{3}{2}}\ga)\int_{M'}\phi^2\,\d V_{g'}.
\end{align*}
Here between the third and fourth lines we have used \eq{Aeq5}, between 
the fifth and sixth lines we have used \eq{Aeq1}, and between the last two 
we have used the $L^2$- orthogonality of $v_1$ and $v_2$ and the trivial
inequality $(b-2\ga)[1+\ga/2(a+b+2\ga)]+a\ga/2(a+b+2\ga)\le b-3\ga/2$. 
This proves \eq{Aeq6}.

To prove \eq{Aeq7}, observe that by \eq{Aeq5}, orthogonal projection 
from $\ti E'$ to $E'$ is injective, and as they have the same (finite) 
dimension, it must also be surjective. Let $\ti v_2\in(\ti E')^\perp$. 
Then $\ti v_2=v_1+v_2$ with $v_1\in E'$ and $v_2\in (E')^\perp$. By this 
surjectivity, there exists $\ti v_1\in\ti E'$ such that $\ti v_1=v_1+v_3$ 
with $v_3\in(E')^\perp$, that is, the $E'$-component of $\ti v_1$ is
$v_1$, the same as that of $\ti v_2$. But $\ti v_1$ and $\ti v_2$ are
orthogonal in $L^2_1(M')$, so taking their inner product gives that
$\snm{v_1}{2}{1}^2=-\an{v_2,v_3}\le\snm{v_3}{2}{1}\snm{v_2}{2}{1}$.

As $\ti v_1=v_1+v_3\in\ti E'$, we may square this inequality,
substitute in for $\snm{v_3}{2}{1}^2$ using \eq{Aeq5}, and divide
through by $\snm{v_1}{2}{1}^2$. The result is that
\begin{equation*}
\snm{v_1}{2}{1}^2\le\frac{\ga}{2(a+b+2\ga)}\snm{v_2}{2}{1}^2,
\quad\text{for $\ti v_2=v_1+v_2 \in(\ti E')^\perp$,}
\end{equation*}
which is the analogue of \eq{Aeq5} for $(\ti E')^\perp$ instead of 
$\ti E'$. This is the ingredient needed to prove \eq{Aeq7} by the 
method used above for \eq{Aeq6}, and the remainder of the proof 
will be left to the reader. 
\end{proof}

For the case of \S\ref{s22}, a subspace $\ti E''$ of functions on $M''$ is 
created in the same way, and Lemma \ref{lA12} applies to this space too. We 
now define the spaces $E_t$. In \S\ref{s21}, let $E_t$ be the 
space of functions that are equal to some function in $\ti E'$ on the 
subset of $M$ identified with $M'\sm\Phi'[\ovB_t(0)]$, and are zero
outside this subset. In \S\ref{s22}, let $E_t$ be the direct
sum of this space of functions, and the corresponding space made
from~$\ti E''$.

For small $t$ the functions in $E_t$ are $C^\iy$, and on their support
$g_t$ is equal to $g'$ (or $g''$ in case \S\ref{s22}). Thus \eq{Aeq6}
applies to functions in~$E_t$.

\begin{lem} For small enough\/ $t$, the subspace
$E_t$ of\/ $L^2_1(M)$ satisfies
\e
\text{if\/ $\phi\in E_t$, then\ }
\int_Ma\ms{\na\phi}\,\d V_{g_t} \le 
\int_M(b-{\ts\frac{3}{2}}\ga)\phi^2\,\d V_{g_t}.
\label{Aeq8}
\e
A fortiori, it satisfies the inequality \eq{Aeq3} of Lemma~\ref{lA11}.
\label{lA13}
\end{lem}

\begin{proof} In \S\ref{s21} this follows immediately from \eq{Aeq6},
as the $g_t$ and $g'$ agree upon the support of the functions of $E_t$.
In \S\ref{s22} $\phi$ is the sum of elements of $\ti E'$ and $\ti E''$;
both sides of \eq{Aeq8} split into two terms, each involving one function.
So \eq{Aeq8} is the sum of two inequalities, which follow immediately from
\eq{Aeq6} as before, and from the counterpart of \eq{Aeq6} applying
to~$\ti E''$.
\end{proof}

The previous lemma showed that the space of functions $E_t$ upon $M$ 
satisfies inequality \eq{Aeq3} of Lemma \ref{lA11}. In the next proposition,
proved in \S\ref{A2}, we show that the inequality \eq{Aeq4} is satisfied too.

\begin{prop} For small\/ $t$, inequality \eq{Aeq4} of Lemma \ref{lA11} holds.
\label{pA14}
\end{prop}

The proof of this proposition is the subject of \S\ref{A2}. Suppose for the 
moment that the proposition holds. Then a space of functions $E_t$ upon $M$ 
has been constructed, satisfying inequality \eq{Aeq3} by Lemma \ref{lA13}, 
and inequality \eq{Aeq4} by Proposition \ref{pA14}. So by Lemma \ref{lA11}, 
the proof of Theorem \ref{t331} is finished. 
\end{proof}

\subsection{Eigenvalues of $a\De$ larger than $b$}
\label{A2}

Next we shall prove Proposition \ref{pA14}, which completes the proof of 
Theorem \ref{t331}. The idea of the proof is as follows. Given $\phi\in
(E_t)^\perp$, we want to show that its `average eigenvalue' of $a\De$
is at least $b+\ga$. We do this in different ways depending on whether
$\phi$ is concentrated in $M'$, or $M''$, or the `neck' in between.

\begin{result}{Proposition \ref{pA14}} Let\/ $M$, $g_t$ and\/ $E_t$
be as in \S\ref{A1}. Then for small\/ $t$, if\/ $\phi\in(E_t)^\perp
\subset L^2_1(M)$, then
\e
\int_Ma\ms{\na\phi}\,\d V_{g_t}\ge\int_M(b+\ga)\phi^2\,\d V_{g_t}.
\label{Aeq9}
\e
\end{result}

\begin{proof} For simplicity we shall prove the proposition for 
the metrics of \S\ref{s21} only, and the modifications for \S\ref{s22}
will be left to the reader. We will start from \eq{Aeq7} of Lemma
\ref{lA12}. The constants in \eq{Aeq9} and \eq{Aeq7} are different
-- the first has $b+\ga$, the second $b+3\ga/2$. Choose constants
$b_0>b_1>\cdots>b_5$ with $b_0=b+3\ga/2$ and $b_5=b+\ga$. These
will be used to contain five error terms.

Shortly we shall choose $r_1,r_2,r_3$ with $0<r_1<r_2<r_3$. Define
three compact Riemannian submanifolds with boundary $R_t\subset S_t
\subset T_t$ in $(M,g_t)$ to be the subsets of $M$ coming from subsets
$R,S$ and $T$ of $M''$ respectively, where
\begin{gather*}
R=\hat M\sm\Xi''\bigl[\R^n\sm\ovB_{r_1}(0)\bigr],\quad
S=\hat M\sm\Xi''\bigl[\R^n\sm\ovB_{r_2}(0)\bigr]\\
\text{and}\quad T=\hat M\sm\Xi''\bigl[\R^n\sm\ovB_{r_3}(0)\bigr].
\end{gather*}

When $t$ is small, $R_t,S_t$ and $T_t$ lie in the region of $M$ in
which the function $\be_2$, used in \S\ref{s21} to define $g_t$, is
1. Then $R_t,S_t,T_t$ are homothetic to $R,S,T$ respectively,
by a homothety multiplying their metrics by~$t^{12}$.

The idea is this. A diffeomorphism $\Psi_t'$ from $M'\sm\{m'\}$ onto 
$M\sm R_t$ will be constructed, which will be the identity outside
$T_t$. Using $\Psi_t'$ any function in $L^2_1(M)$ defines a function
in $L^2_1(M')$. Applying \eq{Aeq7} of Lemma \ref{lA12} therefore
induces an inequality upon functions in $L^2_1(M)$. We will be able
to show that for functions that are not too large in $S_t$, this
inequality implies \eq{Aeq9} as we require. Then only the case of
functions that are large in $S_t$ remains.

Suppose, for the moment, that $r_1,r_2,r_3$ are fixed with $r_1<r_2<r_3$.
For $R_t$ to be well defined, $r_1$ must satisfy $r_1>\de^{-4}$. For
$T_t$ to be well defined, $t$ must be sufficiently small that $t^6r_3<\de$.
We also suppose that $t$ is small enough that the functions of $E_t$
vanish on~$T_t$.

Let $\Psi'_t:M'\sm\{m'\}\ra M$ be the identity outside
$\Phi'\bigl[B_{t^6r_3}(0)\bigr]$ in $M'$, and on
$\Phi'\bigl[B_{t^6r_3}(0)\bigr]$ in $M'$ define $\Psi'_t$ by
\begin{equation*}
\Psi'_t\left(\Phi'(v)\right)=
\Phi'\left(\frac{t^6r_1v}{\md{v}}+\frac{(r_3-r_1)v}{r_3}\right).
\end{equation*}
Let $\phi\in L^2_1(M)$, and define $\phi'=(\Psi'_t)^*(\phi)$. 
Then $\phi'\in L^2_1(M')$, as we shall see.

An easy calculation shows that
\e
b_0\int_{M'}(\phi')^2\,\d V_{g'}=
b_0\int_M\phi^2\cdot F_t\,\d V_{g_t},
\label{Aeq10}
\e
where $F_t$ is a function on $M$ that is 1 on that part of $M$
coming from $M'\sm\Phi'\bigl[B_{t^6r_3}(0)\bigr]$, is zero on that 
part of $M$ not coming from $M'\sm\Phi'\bigl[B_{t^6r_1}(0)\bigr]$, 
and in the intermediate annulus is given by
\begin{equation*}
F_t\left((\phi')^{-1}(v)\right)=
\frac{(\md{v}-t^6r_1)^{n-1}r_3^n}{\md{v}^{n-1}(r_3-r_1)^n}\cdot
\psi'(v)^p\psi_t(v)^{-p}.
\end{equation*}

Similarly, we may easily show that
\e
a\int_{M'}\ms{\na\phi'}\,\d V_{g'}\le
a\int_M\ms{\na\phi}\cdot G_t\,\d V_{g_t},
\label{Aeq11}
\e
where $G_t$ is a function on $M$ that is 1 on that part of $M$ coming 
from $M'\sm\Phi'\bigl[B_{t^6r_3}(0)\bigr]$, is zero on that part of $M$ 
not coming from $M'\sm\Phi'\bigl[B_{t^6r_1}(0)\bigr]$, and in the 
intermediate annulus is given by
\begin{equation*}
G_t\left((\phi')^{-1}(v)\right)\!=\!\max\left(
\frac{(\md{v}\!-\!t^6r_1)^{n\!-\!1}r_3^{n\!+\!2}}{
\md{v}^{n\!-\!1}(r_3\!-\!r_1)^{n\!+\!2}},
\frac{(\md{v}\!-\!t^6r_1)^{n\!-\!3}r_3^{n\!-\!2}}{
\md{v}^{n\!-\!3}(r_3\!-\!r_1)^{n\!-\!2}}
\right)\cdot\psi'(v)\psi_t(v)^{-1}.
\end{equation*}
Here, the first term in the $\max(\dots)$ is the multiplier for the radial 
component of $\na\phi$, and the second term is the multiplier for the 
nonradial components. As $F_t,G_t$ are bounded, we see from \eq{Aeq10} and
\eq{Aeq11} that $\phi'\in L^2_1(M')$, as was stated above.

Suppose now that $\phi\in (E_t)^\perp\subset L^2_1(M)$. For small enough
$t$ this implies that $\phi'\in(\ti E')^\perp$, and so \eq{Aeq7} 
applies by Lemma \ref{lA12}. Combining this with \eq{Aeq10} and \eq{Aeq11}
gives that
\e
a\int_M\ms{\na\phi}\cdot G_t\,\d V_{g_t}\ge
b_0\int_M\phi^2\cdot F_t\,\d V_{g_t}.
\label{Aeq12}
\e

Now by the definition of $\psi_t$, $\psi'(v)\psi_t(v)^{-1}$ approaches
1 as $t\ra 0$. In fact it may be shown that
\begin{equation*}
\bmd{\psi'(v)\psi_t(v)^{-1}-1}\le C_0t^{6(n-2)}\md{v}^{-(n-2)}
\quad\text{when $t^6\le \md{v}\le t^{6-2/(n-2)}$,}
\end{equation*}
for some constant $C_0$. For $t$ small enough this certainly holds
in the region $t^6r_1\le\md{v}\le t^6r_3$, and in this region we have
$\bmd{\psi'(v)\psi_t(v)^{-1}-1}\le C_0^{\phantom{(}}r_1^{-(n-2)}$.

Choose $r_1$ greater than $\de^{-4}$, and large enough that
$b_1(1+C_0r_1^{-(n-2)})^p\le b_0$ and $b_2\le b_1(1-C_0r_1^{-(n-2)})$.
Then for small $t$, the $\psi'\psi_t^{-1}$ terms in $F_t$ and $G_t$
can be absorbed by putting $b_2$ in place of $b_0$. Next, $r_2$ is
defined uniquely in terms of $r_3$ to satisfy $r_1<r_2<r_3$ and
$b_3(r_2-r_1)^{n-1}r_3^nr_2^{1-n}(r_3-r_1)^{-n}=b_4$. Then
$b_3(\md{v}-t^6r_1)^{n-1}r_3^n\md{v}^{1-n}(r_3-r_1)^{-n}\ge b_4$
when $t^6r_2\le\md{v}\le t^6r_3$. This is to bound $F_t$ below on
the region~$\md{v}\ge t^6r_2$.

Finally, we define $r_3$. Choose $r_3>r_1$ sufficiently large that
two conditions hold. The first is that
$b_3\cdot\max\bigl((\md{v}-t^6r_1)^{n-1}r_3^{n+2}\md{v}^{1-n}
(r_3-r_1)^{-(n+2)},(\md{v}-t^6r_1)^{n-3}r_3^{n-2}\md{v}^{3-n}(r_3-r_1)^{2-n}
\bigr)\le b_2$ when $t^6r_1\le\md{v}\le t^6r_3$; combining this with one of
the inequalities used to define $r_1$ shows that $b_3G_t\le b_1$.
The second condition is that
\e
\frac{\vol(S)}{\vol(T)-\vol(S)}\le\frac{b_4-b_5}{4b_5}.
\label{Aeq13}
\e
This condition will be needed later.

The last two definitions are circular, as $r_2$ is defined in terms
of $r_3$, and vice versa because $S_t$ depends on $r_2$. However,
manipulating the definition of $r_2$ reveals that however large $r_3$
is, $r_2$ must satisfy $r_2\le r_1(1-b_3^{1/(n-1)}b_2^{-1/(n-1)})^{-1}$,
and so $\vol(S)$ is bounded in terms of $r_1$, whereas $\vol(T)$ can
grow arbitrarily large. Therefore \eq{Aeq13} does hold for $r_3$
sufficiently large.

The above estimates show that $G_t\le b_1/b_3$ and $F_t\ge b_1b_4/b_0b_3$ 
on $M\sm S_t$ for small $t$. Substituting these into \eq{Aeq12} gives that
when $t$ is sufficiently small,
\e
a\int_M\ms{\na\phi}\,\d V_{g_t}
\ge b_4\int_{M\sm S_t}\phi^2\,\d V_{g_t}.
\label{Aeq14}
\e

Suppose that $\int_{S_t}\phi^2\d V_{g_t}\!\le\!(b_4/b_5-1)\cdot
\int_{M\sm S_t}\phi^2\d V_{g_t}$. Then $b_4\int_{M\sm S_t}\phi^2\d V_{g_t}
\!\ge\!b_5\int_M\phi^2\d V_{g_t}$, and from \eq{Aeq14} we see that 
\eq{Aeq9} holds for $\phi$, which is what we have to prove. Therefore 
it remains only to deal with the case that $\int_{S_t}\phi^2\d V_{g_t}> 
(b_4/b_5-1)\cdot\int_{M\sm S_t}\phi^2\d V_{g_t}$.

Suppose that this inequality holds. The basic idea of the remainder of the
proof is that when $t$ is small, the volume of $S_t$ is also small, and this
forces $\phi$ to be large on $S_t$ compared to its average value elsewhere.
Therefore $\phi$ must change substantially in the neighbourhood $T_t$ of
$S_t$, and this forces $\na\phi$ to be large in $T_t$.

Restrict $t$ further, to be small enough that $t^6r_3\le t^2$. Then $T_t$ is
contained in the region of gluing in which $\be_2=1$. So the pair
$(S_t,T_t)$ is homothetic to a pair $(S,T)$ of compact manifolds with
$C^\iy$ boundaries and with $S$ contained in the interior of $T$; the
metrics on $(S_t,T_t)$ are the metrics on $(S,T)$ multiplied by $t^{12}$. For
these $S,T$ the following lemma holds.

\begin{lem} Let\/ $S,T$ be compact, connected Riemannian 
$n$-manifolds with smooth boundaries, such that\/ $S\subset T$
but\/ $S\ne T$. Then there exists\/ $C_1>0$ such that
for all\/ $\phi\in L^2_1(T)$,
\e
\left(\frac{\int_S\phi^2\,\d V_g}{\vol(S)}\right)^{1/2}-
\left(\frac{\int_{T\sm S}\phi^2\,\d V_g}{\vol(T)-\vol(S)}\right)^{1/2}
\le C_1\left(\int_T\ms{\na\phi}\,\d V_g\right)^{1/2}.
\label{Aeq15}
\e
\label{lA21}
\end{lem}

\begin{proof} We begin by quoting a theorem on the existence of
solutions of the equation $\De u=f$ on a manifold with smooth boundary.

\begin{thm} Suppose that\/ $T$ is a compact manifold with smooth
boundary, and that\/ $f\in L^2(T)$ and satisfies $\int_Tf\d V_g=0$.
Then there exists $\xi\in L^2_2(T)$, unique up to the addition of a
constant, such that\/ $\De\xi=f$, and in addition ${\bf n}\cdot\na\xi$
vanishes at the boundary, where ${\bf n}$ is the unit outward normal
to the boundary.
\end{thm}

\begin{proof} This is a partial statement of \cite[Ex.~2, p.~65]{Horm}. 
H\"ormander's example is only stated for $C^\iy$ functions $f$
and $\xi$, but the proof works for $f\in {\cal H}_{\sst (0)}(T)$
and $\xi\in {\cal H}_{\sst (2)}(T)$ in his notation, which are
$L^2(T)$ and $L^2_2(T)$ in ours. 
\end{proof}

Put $f=\vol(S)^{-1}$ in $S$ and $f=(\vol(S)-\vol(T))^{-1}$ in $T\sm S$. 
Then $\int_Tf\d V_g=0$, so by the theorem, there exists a function $\xi\in
L^2_2(M)$ satisfying $\De\xi=f$, and that $\na\xi$ vanishes normal to
the boundary. Because of this vanishing, the boundary term has dropped out 
of the following integration by parts equation:
\begin{equation*}
\int_T\phi\De\xi=-\int_T\bigl(\na\phi\bigr)\cdot
\bigl(\na\xi\bigr)\,\d V_g.
\end{equation*}
Substituting in for $\De\xi$ and using H\"older's inequality gives
\e
\begin{split}
\frac{1}{\vol(S)}\left\vert\int_S\phi\,\d V_g\right\vert
-\frac{1}{\vol(T)-\vol(S)}\left\vert\int_{T\sm S}\phi\,\d V_g\right\vert&\\
\le \left(\int_T\ms{\na\xi}\,\d V_g\right)^{1/2}\cdot
\left(\int_T\ms{\na\phi}\,\d V_g\right)^{1/2}&.
\end{split}
\label{Aeq16}
\e

Now $S$ is connected, so the constants are the only eigenvectors of $\De$
on $S$ with eigenvalue 0 and derivative vanishing normal to the boundary. By
the discreteness of the spectrum of $\De$ on $S$ with these boundary
conditions, there is a positive constant $K_S$ less than or equal to all
the positive eigenvalues. It easily follows that for $\phi\in L^2_1(S)$, 
\e
\left(\frac{\int_S\phi^2\,\d V_g}{\vol(S)}\right)^{1/2}\le
\frac{1}{\vol(S)}\left\vert\int_S\phi\,\d V_g\right\vert+
\left(\frac{\int_S\ms{\na\phi}\,\d V_g}{K_S\cdot\vol(S)}\right)^{1/2}.
\label{Aeq17}
\e
Also, a simple application of H\"older's inequality yields
\e
\frac{1}{\vol(T)-\vol(S)}\left\vert\int_{T\sm S}\phi\,\d V_g\right\vert\le
\left(\frac{\int_{T\sm S}\phi^2\,\d V_g}{\vol(T)-\vol(S)}\right)^{1/2}.
\label{Aeq18}
\e

Adding together \eq{Aeq16}, \eq{Aeq17} and \eq{Aeq18} gives \eq{Aeq15}, 
as we want, with constant $C_1=\bigl(\int_T\ms{\na\xi}\,\d V_g
\bigr)^{1/2}+(K_S\cdot\vol(S))^{-1/2}$. 
\end{proof}

The point of this calculation is that because $(S_t,T_t)$ are homothetic to
$(S,T)$ by the constant factor $t^{12}$, Lemma \ref{lA21} implies that for all
$\phi\in L^2_1(T_t)$,
\begin{equation*}
\left(\frac{\int_{S_t}\phi^2\,\d V_{g_t}}{\vol(S)}\right)^{1/2}-\left(
\frac{\int_{T_t\sm S_t}\phi^2\,\d V_{g_t}}{\vol(T)-\vol(S)}\right)^{1/2}
\le C_1t^6\left(\int_{T_t}\ms{\na\phi}\,\d V_{g_t}\right)^{1/2}.
\end{equation*}

Now, using \eq{Aeq13} it follows that
\e
\begin{split}
\left(\int_{S_t}\phi^2\,\d V_{g_t}\right)^{1/2}-\left(
\frac{b_4-b_5}{4b_5}\int_{M\sm S_t}\phi^2\,\d V_{g_t}\right)^{1/2}&\\
\le \vol(S)^{1/2}C_1t^6
\left(\int_M\ms{\na\phi}\,\d V_{g_t}\right)^{1/2}&.
\end{split}
\label{Aeq19}
\e

But because $\int_{S_t}\phi^2\d V_{g_t}>(b_4/b_5-1)\cdot
\int_{M\sm S_t}\phi^2\d V_{g_t}$, substituting this into
\eq{Aeq19}, squaring and manipulating gives that
\begin{equation*}
b_5\int_M\phi^2\,\d V_{g_t}<
\frac{b_4b_5}{b_4-b_5}\cdot\int_{S_t}\phi^2\,\d V_{g_t}
<\frac{4b_4b_5t^{12}C_1^2\vol(S)}{b_4-b_5}
\cdot\int_M\ms{\na\phi}\,\d V_{g_t}.
\end{equation*}
Therefore, if $t$ is small, then inequality \eq{Aeq9} holds. This
completes the proof of Proposition~\ref{pA14}.
\end{proof}

\subsection{The spectrum of $a\De$ in the zero scalar curvature cases}
\label{A3}

Now we prove Proposition \ref{p433} and Theorems \ref{t432} and \ref{t434}.
We shall start with a preliminary version of Theorems \ref{t432} and
\ref{t434}. Suppose $\ga>0$ such that $a\De$ has no eigenvalues in
$(0,2\ga)$ on $(M',g')$ in \S\ref{s41}, and on both $(M',g')$ and
$(M'',g'')$ in \S\ref{s42}. Then
\begin{equation*}
\text{if $\phi\in L^2_1(M')$ and $\int_{M'}\phi
\,\d V_{g'}=0$, then\ }
\int_{M'}a\ms{\na\phi}\,\d V_{g'} \ge 
\int_{M'}2\ga\phi^2\,\d V_{g'},
\end{equation*}
and the same for $M''$ in \S\ref{s42}. This is an analogue of 
\eq{Aeq2} of \S\ref{A1}, and the analogues of $E'$ and $E''$ are
the spaces of constant functions on $M'$ and~$M''$.

Define spaces of functions $E_t$ on $M$ as in \S\ref{A1}. This involves
choosing a function $\si'$ on $M'$ (done just before Lemma \ref{lA12}) 
that is 1 outside a small ball around $m'$, and a similar 
function $\si''$ on $M''$. For the proof of Proposition \ref{p433} 
later, we must choose $\si',\si''$ to vary with $t$ rather 
than being fixed as they were before. Let $\si',\si''$ be
functions on $M',M''$ defined in the same way as $\be_1,\be_2$
in \S\ref{s41}, but which are 0 when $\md{v}\le t^{2k/(n+2)}$ and
1 when $\md{v}\ge t^{(n-2)/(n+1)}$ and outside~$B',B''$.

This makes sense for small $t$, as the definition of $k$ in
\S\ref{s41} gives $2k/(n+2)>(n-2)/(n+1)$ so that $t^{2k/(n+2)}<
t^{(n-2)/(n+1)}<\de$ for small $t$. Note also that the choice of
the inner radius $t^{2k/(n+2)}$ means (from \S\ref{s41}) that
$g_t$ is identified with $g'$ and $g''$ on the support of
$\si'$ and $\si''$ respectively, so that the functions $E_t$ are
supported in the parts of $M$ where $g_t$ equals $g'$ or $g''$.
Here is a first approximation to Theorems \ref{t432} and~\ref{t434}:

\begin{lem} Let\/ $\{g_t:t\in (0,\de)\}$ be one of the families of 
metrics defined in \S\ref{s41} or \S\ref{s42}. Then for small\/~$t$,
\begin{equation*}
\text{if\/ $\phi\in (E_t)^\perp\subset L^2_1(M)$, then\ }
\int_Ma\ms{\na\phi}\,\d V_{g_t} \ge 
\ga\int_M\phi^2\,\d V_{g_t}.
\end{equation*}
Here the orthogonal space is taken with respect to the $L^2$ inner product.
\label{lA31}
\end{lem}

\begin{proof} This is proved just as is Proposition \ref{pA14}, except that 
we use the inner product in $L^2(M)$ rather than that in $L^2_1(M)$, and
instead of choosing a series of constants interpolating between $b+2\ga$
and $b+\ga$, we choose constants interpolating between $2\ga$ and $\ga$,
and some simple changes must be made to the proof as the powers of $t$
used in defining the metrics of \S\ref{s41} and \S\ref{s42} are different
to those used in \S\ref{s21} and~\S\ref{s22}.
\end{proof}

Now as $E_t$ is modelled on $E'$ in \S\ref{s41} and on
$E'\oplus E''$ in \S\ref{s42}, it is close to the constant 
functions in \S\ref{s41}, and to functions taking one constant
value on the $M'$ part of $M$ and another constant value on the
$M''$ part in~\S\ref{s42}.

These spaces $E_t$ are not quite good enough for our purposes,
for we shall need spaces that contain the constants. We therefore
produce modified spaces $\ti E_t$, and prove a similar lemma for them.
In \S\ref{s41}, let $\ti E_t=\an{1}$. In \S\ref{s42}, let $e\in E_t$
be the unique element that is nonnegative on the part of $M$ coming
from $M'$ and satisfies $\int_Me\,\d V_{g_t}=0$ and $\int_Me^2\d
V_{g_t}=2\vol(M')$, and let~$\ti E_t=\an{1,e}$.

\begin{lem} Let\/ $\{g_t:t\in (0,\de)\}$ be one of the families of 
metrics defined in \S\ref{s41} or \S\ref{s42}. Then for small\/~$t$,
\begin{equation*}
\text{if\/ $\phi\in(\ti E_t)^\perp\subset L^2_1(M)$, then\ }
\int_Ma\ms{\na\phi}\,\d V_{g_t} \ge 
\ga\int_M\phi^2\,\d V_{g_t}.
\end{equation*}
Here the orthogonal space is taken with respect to the $L^2$ inner product.
\label{lA32}
\end{lem}

\begin{proof} Let $\xi$ be the unique element of $E_t$ with
$\int_M\xi\d V_{g_t}=1$, and $\int_Me\xi\d V_{g_t}=0$ in \S\ref{s42}.
If $\phi\in(\ti E_t)^\perp$, then $\phi-\an{\phi,\xi}\in
E_t^\perp$, where $\an{\,,\,}$ is the inner product of $L^2(M)$.
So by Lemma \ref{lA31},
\begin{align*}
\int_Ma\ms{\na\phi}\,\d V_{g_t} &\ge \ga
\int_M\bigl(\phi-\an{\phi,\xi}\bigr)^2\,\d V_{g_t}\\
&=\ga\int_M\bigl(\phi^2+\an{\phi,\xi}^2\bigr)\,\d V_{g_t}\\
&\ge\ga\int_M\phi^2\,\d V_{g_t}.
\end{align*}
Here between the first and second lines we have used the fact that
$\phi\in(\ti E_t)^\perp$, and thus $\int_M\phi\,\d V_{g_t}=0$, as 
the constants lie in $\ti E_t$. 
\end{proof}

As $\ti E_t$ is just the constant functions in \S\ref{s21}, we
immediately deduce 

\begin{cor} Theorem \ref{t432} is true.
\label{cA33}
\end{cor}

We will now construct the eigenvector $\be_t$ of Proposition \ref{p433}, using a
sequence method. The proposition is stated again here.

\begin{result}{Proposition \ref{p433}} Let\/ $\{g_t:t\in (0,\de)\}$
be the family of metrics defined on $M=M'\# M''$ in \S\ref{s42}. Then
for small\/ $t$, there exists $\la_t>0$ and\/ $\be_t\in C^\iy(M)$ such
that\/ $a\De\be_t=\la_t\be_t$ on $(M,g_t)$. Here $\la_t=O(t^{n-2})$, and
\e
\be_t=\begin{cases}
\phantom{-}1+O(t^{n-2}) & \text{on $M'\sm B'$}\\
\phantom{-}1+O(t^{n-2}\md{v}^{2-n})
&\text{on $\{v:t\le\md{v}<\de\}\subset B'$} \\
-1+O(t^{n-2}) & \text{on $M''\sm B''$} \\
-1+O(t^{n-2}\md{v}^{2-n}) 
& \text{on $\{v:t\le\md{v}<\de\}\subset B''$,}
\end{cases}
\label{Aeq20}
\e
identifying subsets of\/ $M',M''$ with subsets of\/~$M$.
\end{result}

\begin{proof} Let $y$ be the unique element of $C^\iy(M)$
satisfying $\int_My\,\d V_{g_t}=0$ and $a\De y=e$. To make $\be_t=e+w$ 
and $a\De\be_t=\la_t\be_t$, we must find $w$ and $\la_t$ such that $a\De
w=\la_t(e+w)-a\De e$. Define inductively a sequence of real numbers
$\{\la_i\}_{i=0}^\iy$ and a sequence $\{w_i\}_{i=0}^\iy$ of
elements of $C^\iy(M)$ beginning with $\la_0=w_0=0$. Let
\e
\la_i=2\vol(M')\left(\int_My(e+w_{i-1})\,\d V_{g_t}\right)^{-1},
\label{Aeq21}
\e
and let $w_i$ be the unique element of $C^\iy(M)$ satisfying
$\int_Mw_i\,\d V_{g_t}=0$ and 
\e
a\De w_i=\la_i(e+w_{i-1})-a\De e.
\label{Aeq22}
\e
Note that as $\int_Me\,\d V_{g_t}=\int_Mw_{i-1}\,\d V_{g_t}=0$, the right
hand side has integral zero over $M$, and so $w_i$ exists. Thus the 
sequences $\{\la_i\}_{i=0}^\iy$ and $\{w_i\}_{i=0}^\iy$ are 
well-defined, provided only that the integral on the right hand side of 
\eq{Aeq21} is nonzero; we will prove later that the integral is bounded
below by a positive constant.

If both sequences converge to $\la_t$ and $w$ respectively, say, then 
\eq{Aeq22} implies that $a\De w=\la_t(e+w)-a\De e$, so that $\be_t=e+w$ is
an eigenvector of $a\De$ associated to $\la_t$. The r\^ole of \eq{Aeq21}
is as follows: multiply \eq{Aeq22} by $y$ and integrate over $M$.
Integrating by parts gives
\begin{equation*}
\int_Mw_ia\De y\,\d V_{g_t}=\la_i\int_My(e+w_{i-1})\d V_{g_t}
-\int_Mea\De y\,\d V_{g_t}.
\end{equation*}

Substituting $e$ for $a\De y$ and recalling that
$\int_Me^2\d V_{g_t}=2\vol(M')$, we get
\begin{equation*}
\int_Mew_i\,\d V_{g_t}=\la_i\int_My(e+w_{i-1})\,\d V_{g_t}
-2\vol(M')=0,
\end{equation*}
so that $\int_Mew_i\d V_{g_t}=0$, and if $w$ is the limit of the
sequence then $\int_Mew\,\d V_{g_t}\!=\nobreak\!\nobreak 0$.

Now as $\int_Mw_i\d V_{g_t}=\int_Mew_i\d V_{g_t}=0$, 
$w_i\in(\ti E_t)^\perp$ and Lemma \ref{lA32} gives
\e
a\lnm{\na w_i}{2}^2\ge\ga\lnm{w_i}{2}^2.
\label{Aeq23}
\e
Multiplying \eq{Aeq22} by $w_i$ and integrating over $M$ by parts gives
\begin{align*}
a\lnm{\na w_i}{2}^2&=\int_Ma\ms{\na w_i}\,\d V_{g_t}\\
&=\la_i\int_Mw_i(e+w_{i-1})\,\d V_{g_t}-a\int_M\na w_i\cdot\na 
e\,\d V_{g_t}\\
&\le \md{\la_i}\lnm{w_i}{2}\bigl(\lnm{e}{2}+\lnm{w_{i-1}}{2}\bigr)
+a\lnm{\na w_i}{2}\lnm{\na e}{2}\\
&\le\md{\la_i}a^{1/2}\ga^{-1/2} \lnm{\na w_i}{2}
\bigl(\lnm{e}{2}+\lnm{w_{i-1}}{2}\bigr)+a\lnm{\na w_i}{2}\lnm{\na
e}{2},
\end{align*}
applying \eq{Aeq23} between the second and third lines, and H\"older's
inequality. Dividing by $a^{1/2}\ga^{-1/2}\lnm{\na w_i}{2}$ and using
\eq{Aeq23} on the left hand side gives
\e
\ga\lnm{w_i}{2}\le
\md{\la_i}\bigl(\lnm{e}{2}+\lnm{w_{i-1}}{2}\bigr)+D_0,
\label{Aeq24}
\e
where $D_0=(a\ga)^{1/2}\lnm{\na e}{2}$.

Define $D_1=(2\vol(M'))^{-1}\int_Mye\,\d V_{g_t}$. Then $y=D_1e+z$, where
$\int_Mz\,\d V_{g_t}=\int_Mze\,\d V_{g_t}=0$. So $z\in(\ti E_t)^\perp$,
and by Lemma~\ref{lA32}
\e
a\lnm{\na z}{2}^2\ge\ga\lnm{z}{2}^2.
\label{Aeq25}
\e

As $\int_Mze\,\d V_{g_t}=0$ and $e=a\De y=aD_1\De e+a\De z$ we have
$\int_Mz\bigl(D_1\De e+\De z\bigr)\d V_{g_t}=0$, so
\begin{equation*}
\lnm{\na z}{2}^2=-D_1\int_M\na z\cdot\na e\,\d V_{g_t}
\le D_1\lnm{\na z}{2}\lnm{\na e}{2}.
\end{equation*}
Multiplying by $(a\ga)^{1/2}\lnm{\na z}{2}^{-1}$ and substituting 
\eq{Aeq25} into the l.h.s.\ then gives $\ga\lnm{z}{2}\le D_0D_1$. Since
$y=D_1e+z$ and $\int_Mew_{i-1}\d V_{g_t}=0$, from \eq{Aeq21} we find
\e
\begin{split}
\md{\la_i}^{-1}=D_1+\frac{1}{2\vol(M')}\int_Mzw_{i-1}\,\d V_{g_t}
&\ge D_1-\frac{\lnm{z}{2}\lnm{w_{i-1}}{2}}{2\vol(M')}\\
&\ge D_1-\frac{D_0D_1\lnm{w_{i-1}}{2}}{2\vol(M')\ga}\,.
\end{split}
\label{Aeq26}
\e

Now \eq{Aeq24} and \eq{Aeq26} are what we need to prove that the
sequences $\{\la_i\}_{i=0}^\iy$ and $\{w_i\}_{i=0}^\iy$ are well-defined
and convergent, provided $D_0$ is sufficiently small and $D_1$
sufficiently large. It can be shown that if $2D_0^2\le\ga^2\vol(M')$
and $D_1$ is large enough, then the two sequences converge to $\la_t$
and $w$ respectively satisfying
\begin{gather*}
a\De w=\la_t(e+w)-a\De e,
\quad\text{where}\quad \lnm{w}{2}\le 2D_0/\ga\quad\text{and}\\
\md{\la_t}^{-1}\ge D_1-\frac{D_0^2D_1}{\ga^2\vol(M')}\ge\ha D_1,
\quad\text{so that}\quad \md{\la_t}\le 2D_1^{-1}.
\end{gather*}
The proof uses the same sort of reasoning as Lemma \ref{l312},
and will be left to the reader.

Let us now look more closely at $D_0$ and $D_1$. Firstly,
$D_0=(a\ga)^{1/2}\lnm{\na e}{2}$, and $e$ is defined by 
$\si'$ and $\si''$ chosen just before Lemma \ref{lA31}.
In fact, $e=c'\si'-c''\si''$, where $c'$ and $c''$ are
close to 1, as $\vol(M')=\vol(M'')$. But $\si',\si''$
are defined in the same way as $\be_1,\be_2$ of \S\ref{s21},
and an estimate analogous to \eq{s2eq4} of \S\ref{s24} applies
to them, from which we find that $D_0=O(t^{(n-2)^2/(n+1)})$ for
small $t$. Because $\lnm{w}{2}\le 2D_0/\ga$, this gives
that~$\lnm{w}{2}=O(t^{(n-2)^2/(n+1)})$.

Secondly, we need to estimate $D_1$. Let $\xi'$ be the Green's
function of $a\De$ at $m'$ on $M'$, satisfying $a\De\xi'=\de_{m'}-
1/\vol(M')$ in the sense of distributions. Then $\xi'$ has a pole
of the form $D_2\md{v}^{2-n}+O''(\md{v}^{1-n})$ at $m'$, where
$D_2=(n-2)^{-1}\om_{n-1}^{-1}$. Similarly, we may define the
Green's function $\xi''$ of $a\De$ at $m''$ on~$M''$.

The application of this is in modelling the function $y$ on $M$. We may 
view $M$, approximately, as being made up of the unions of $M'$ and
$M''$, each with a small ball of radius $t$ cut out. To get a
function $\xi$ on $M$ with $\De\xi$ close to $\vol(M')^{-1}$ on the
part coming from $M'$, and close to $-\vol(M'')^{-1}$ on
the part coming from $M''$, we try $\xi$ equal approximately to
$d'-\xi'$ on $M'$ and $\xi''-d''$ on the part coming from
$M''$, for constants $d'$ and $d''$. To join
these two functions together on the neck we must have
$d'+d''=2D_2t^{2-n}+O(t^{1-n})$, and for $\int_M\xi
\d V_{g_t}=0$ we must have $d'\vol(M')=d''\vol(M'')+O(t^{1-n})$. 
As $\vol(M')=\vol(M'')$ this gives~$d'=D_2t^{2-n}+O(t^{1-n})=d''$.

But $e$ is approximately equal to 1 on $M'$ and to $-1$ on
$M''$, so that $e\sim \vol(M')a\De\xi$. Therefore
$y\sim \vol(M')\xi$, and $D_1=(2\vol(M'))^{-1}\int_Mye\,\d V_{g_t}
\sim d'\vol(M')$. So finally, we conclude that
$D_1=D_2\vol(M')t^{2-n}+O(t^{1-n})$. This validates the claim
that $D_1$ is large for small $t$, that was used earlier to 
ensure convergence.

Taking the limit over $i$ in \eq{Aeq21}, we find that
\begin{equation*}
\la_t^{-1}=D_1+\frac{1}{2\vol(M')}\int_Mzw\,\d V_{g_t},
\end{equation*}
so using estimates on $D_1$, $z$ and $w$ gives that $\la_t=O(t^{n-2})$ 
for small $t$, one of the conclusions of the proposition. Also, if $e$
is a first approximation to $\be_t$, then $\la_t y$ is the second, and
the model of $y$ above gives a model of $\be_t$. So we see that
\eq{Aeq20} holds for $\be_t$, completing the proof.
\end{proof}

Finally, we show that Lemma \ref{lA32} may be modified further, to apply to
functions orthogonal to both 1 and the eigenvector $\be_t$ constructed in the
last proposition:

\begin{lem} Let\/ $\{g_t:t\in (0,\de)\}$ be the family of metrics
defined on $M=M'\# M''$ in \S\ref{s42}. Then for small\/ $t$, if\/
$\phi\in L^2_1(M)$ satisfies $\an{\phi,1}=\an{\phi,\be_t}=0$ in
either the $L^2$ or the $L^2_1$ inner product, then
\e
\int_Ma\ms{\na\phi}\,\d V_{g_t} \ge 
\ga\int_M\phi^2\,\d V_{g_t}.
\label{Aeq27}
\e
Here $\be_t$ is the eigenvector of\/ $a\De$ constructed in
Proposition~\ref{p433}.
\label{lA34}
\end{lem}

\begin{proof} The proof is almost the same as that of Lemma \ref{lA32}. 
Note that as $1,\be_t$ are eigenvectors of $\De$, orthogonality to them 
with respect to the $L^2$ and $L^2_1$ inner products is equivalent, and
so we may suppose that $\an{\,,\,}$ is the inner product of $L^2(M)$.
Let $\xi$ be the unique element of $\ti E_t$ satisfying
$\an{\xi,\be_t}=1$ and $\int_M\xi\d V_{g_t}=0$. If $\phi\in L^2_1(M)$
satisfies $\an{\phi,1}=\an{\phi,\be_t}=0$, then
$\phi-\an{\phi,\xi}\be_t\in\ti E_t^\perp$, taken with respect to
the $L^2$ norm. So by Lemma \ref{lA32},
\begin{equation*}
\int_Ma\bms{\na\phi-\an{\phi,\xi}\na\be_t}\,\d V_{g_t}\ge\ga
\int_M\bigl(\phi-\an{\phi,\xi}\be_t\bigr)^2\,\d V_{g_t}.
\end{equation*}
But as $\be_t$ is an eigenvalue of $\De$, it is orthogonal to $\phi$
in both $L^2$ and $L^2_1$ norms, so this equation becomes
\e
\int_Ma\Bigl(\ms{\na\phi}+\an{\phi,\xi}^2\ms{\na\be_t}
\Bigr)\,\d V_{g_t}\ge\ga
\int_M\bigl(\phi^2+\an{\phi,\xi}^2\be_t^2\bigr)\,\d V_{g_t}.
\label{Aeq28}
\e
Now $\la_t<\ga$ for small $t$, so that $\int_Ma\ms{\na\be_t}\d V_{g_t}
\le\ga\int_M\be_t^2\d V_{g_t}$, and subtracting this multiplied by
$\an{\phi,\xi}^2$ from \eq{Aeq28} gives~\eq{Aeq27}. 
\end{proof}

\begin{cor} Theorem \ref{t434} is true.
\label{cA35}
\end{cor}


\begin{thebibliography}{99}

\bibitem{Aubi} T. Aubin, {\it Nonlinear Analysis on Manifolds.
Monge--Amp\`ere Equations}, Grundlehren der math. Wiss. 252, 
Springer-Verlag, 1982.

\bibitem{Bart} R. Bartnik, {\it The Mass of an Asymptotically Flat
Manifold}, Comm. Pure Appl. Math. 39 (1986), 661-693.

\bibitem{GrLa} M. Gromov and H.B. Lawson, {\it The classification of
simply connected manifolds of positive scalar curvature}, Ann. Math.
111 (1980), 423--434.

\bibitem{Horm} L. H\"ormander, {\it Linear Partial Differential Operators},
Grundlehren der math. Wiss. 116, Springer-Verlag, 1963.

\bibitem{Joyc} D.D. Joyce, {\it Hypercomplex and Quaternionic Manifolds,
and Scalar Curvature on Connected Sums}, Oxford D.Phil. thesis, 1992.

\bibitem{Koba} O. Kobayashi, {\it Scalar curvature of a metric with
unit volume}, Math. Ann. 279 (1987), 253--265.

\bibitem{LePa} J.M. Lee and T.H. Parker, {\it The Yamabe Problem}, 
Bull. A.M.S. 17 (1987), 37--91.

\bibitem{MPU} R. Mazzeo, D. Pollack and K. Uhlenbeck, {\it Connected
sum constructions for constant scalar curvature metrics}, Topological
methods in Nonlinear Analysis 6 (1995), 207--233. dg-ga/9511018.

\bibitem{Scho1} R. Schoen, {\it Conformal deformation of a
Riemannian metric to constant scalar curvature}, J. Diff. 
Geom. 20 (1984), 479--495.

\bibitem{Scho2} R. Schoen, {\it Variational theory for the Total
Scalar Curvature Functional for Riemannian Metrics and Related Topics},
pages 120--154 in M. Giaquinta, editor, {\it Topics in the Calculus of
Variations}, Lecture Notes in Mathematics 1365, Springer, 1987.

\bibitem{ScYa} R. Schoen and S.-T. Yau, {\it On the structure of
manifolds with positive scalar curvature}, Manuscripta Math. 28
(1979), 159--183.

\bibitem{Trud} N.S. Trudinger, {\it Remarks concerning the conformal 
deformation of Riemannian structures on compact manifolds}, 
Ann. Scuola Norm. Sup. Pisa 22 (1968), 265--274.

\bibitem{Yama} H. Yamabe, {\it On a deformation of Riemannian structures
on compact manifolds}, Osaka Math. J. 12 (1960), 21--37.

\bibitem{Witt} E. Witten, {\it A new proof of the positive energy
theorem}, Comm. Math. Phys. 80 (1981), 381--402.

\end{thebibliography}
\end{document}